\documentclass{article}
\usepackage{amsmath}
\usepackage{amssymb}
\usepackage{caption}
\usepackage{subcaption}
\usepackage{graphicx}
\usepackage{mathtools}
\usepackage{arydshln}
\usepackage{amsthm}
 \usepackage{relsize}

\usepackage{authblk}
\usepackage{multicol}
\usepackage[final]{changes}

\usepackage{algpseudocode}
\usepackage{algorithm}

\newtheorem{theorem}{Theorem}[section]

\newtheorem{remark}[theorem]{Remark}

\usepackage{xcolor}

\usepackage[utf8]{inputenc}
\usepackage{authblk}

\title{On the efficient preconditioning of the Stokes equations in tight geometries.}

\author[1]{Vladislav Pimanov} 
\author[1]{Oleg Iliev}
\author[2,3]{Ivan Oseledets} 
\author[2]{Ekaterina Muravleva}

\affil[1]{Fraunhofer ITWM, Kaiserslautern, Germany}
\affil[2]{Skolkovo Institute of Technology, Moscow, Russia}
\affil[3]{Artificial Intelligence Research Institute, Moscow, Russia}
\affil[ ]{\textit{Emails:} vladislav.pimanov@itwm.fraunhofer.de, oleg.iliev@itwm.fraunhofer.de, i.oseledets@skoltech.ru, e.muravleva@skoltech.ru}

\makeatletter
\renewcommand\@date{{%
  \vspace{-\baselineskip}%
  \large\centering
  \today
}}
\makeatother

\begin{document}

\maketitle



\begin{abstract}

If the Stokes equations are properly discretized, it is known that the Schur complement matrix is spectrally equivalent to the identity matrix. Moreover, in the case of simple geometries, it is often observed that most of its eigenvalues are equal to one. These facts form the basis for the famous Uzawa algorithm. Despite recent progress in developing efficient iterative methods for solving the Stokes problem, the Uzawa algorithm remains popular in science and engineering, especially when accelerated by Krylov subspace methods.
However, in complex geometries, the Schur complement matrix can become severely ill-conditioned, having a significant portion of non-unit eigenvalues. This makes the established Uzawa preconditioner inefficient.
To explain this behaviour, we examine the Pressure Schur Complement formulation for the staggered finite-difference discretization of the Stokes equations. Firstly, we conjecture that the no-slip boundary conditions are the reason for non-unit eigenvalues of the Schur complement matrix. Secondly, we demonstrate that its condition number increases with increasing the surface-to-volume ratio of the flow domain.

As an alternative to the Uzawa preconditioner, we propose using the diffusive SIMPLE preconditioner for geometries with a large surface-to-volume ratio. We show that the latter is much more fast and robust for such geometries. Furthermore, we show that the usage of the SIMPLE preconditioner leads to more accurate practical computation of the permeability of tight porous media.


\end{abstract}

\section{Introduction}

Let's consider a bounded, open domain $\Omega^f \subset \mathbb{R}^d$, where $d=2,3$. This domain represents the region in which the fluid flow propagates. Within $\Omega^f$, we consider the steady-state Stokes equations, given as follows:
\begin{equation}\label{stokes}
    \begin{split}
        -\Delta \mathbf{u} + \nabla p  &= \mathbf{f} \text{ in } \Omega^f,\\
        -\nabla \cdot \mathbf{u} &= 0 \text{ in } \Omega^f, \\
    \end{split}
\end{equation}
where $\mathbf{u}$ is the fluid velocity, $p$ is the fluid pressure, and $\mathbf{f}$ is the volumetric force driving the flow. 
The Boundary Value Problem (BPV) for the Stokes equations \eqref{stokes} is obtained by additionally imposing  boundary conditions on the boundary $\partial\Omega^f = \Gamma_0$. We consider the no-slip boundary condition on $\Gamma_0$, which is simply the zero Dirichlet condition imposed on the velocity:
\begin{equation}\label{stokes_bc_noslip}
    \mathbf{u} = \mathbf{0} \text{ on } \Gamma_0.
\end{equation}
The discretization of the BVP \eqref{stokes},\eqref{stokes_bc_noslip} leads to a block system of linear equations of the following form:
\begin{equation}\label{coupled_matrix}
\mathbb{A}
    \begin{bmatrix}
    \mathbf{u}_h \\ 
    p_h
    \end{bmatrix}
    =
    \begin{bmatrix}
    \mathbf{f}_h \\ 
     0
    \end{bmatrix}, \quad
    \mathbb{A} = 
    \begin{bmatrix}
    \mathbf{A} & \mathbf{B}^T \\ 
    \mathbf{B} & 0
    \end{bmatrix}.
\end{equation}
Here the vectors $\mathbf{u}_h$, $p_h$, $\mathbf{f}_h$ denote the discretized velocity, pressure, force, respectively, and the subscript $h$ denotes the mesh size parameter. The matrices $\mathbf{A}$ and $\mathbf{B}$ are discrete counterparts of the negative velocity Laplacian operator and the negative divergence operator. The matrix $\mathbf{B}^T$ denotes the discrete pressure gradient operator, which is adjoint to the negative divergence operator under proper discretization. Importantly, the no-slip boundary condition \eqref{stokes_bc_noslip} is incorporated into the matrix $\mathbb{A}$.

In what follows, we consider the Pressure Schur Complement formulation, according to the terminology from \cite{turek1999efficient}, which reduces the coupled system \eqref{coupled_matrix} to an equivalent system for the pressure:
\begin{equation}\label{schur_system}
    S p_h = g_h,
\end{equation}
where $S$ is the Schur complement of the matrix $\mathbb{A}$ and $g_h$ is the right-hand-side in the reduced equation, which are defined as follows:
\begin{equation}\label{schur_rhs}
     S = \mathbf{B}\mathbf{A}^{-1} \mathbf{B}^T, \quad g_h = \mathbf{B}\mathbf{A}^{-1}\mathbf{f}_h.
\end{equation}
Once the pressure is computed, the velocity can be recovered by solving the following system for the velocity:
\begin{equation}
\label{velocity_schur}
    \mathbf{A}\mathbf{u}_h = \mathbf{f}_h - \mathbf{B}^Tp_h.
\end{equation}
In this paper, we examine the impact of the no-slip boundary condition \eqref{stokes_bc_noslip} on the conditioning of the Schur complement system \eqref{schur_system}, particularly when the surface-to-volume ratio of the flow domain $\Omega^f$ is high. The surface-to-volume ratio refers to the ratio between the surface area of the no-slip boundary $\Gamma_0$ and the volume of the flow domain $\Omega^f$.
Note, high surface-to-volume ratio is a distinctive feature for our critical application — namely, computation of permeability of tight porous media. Importantly, our findings are not confined to this specific application.

\begin{figure}
     \centering
     \begin{subfigure}[b]{0.25\textwidth}
         \centering
         \includegraphics[width=\textwidth]{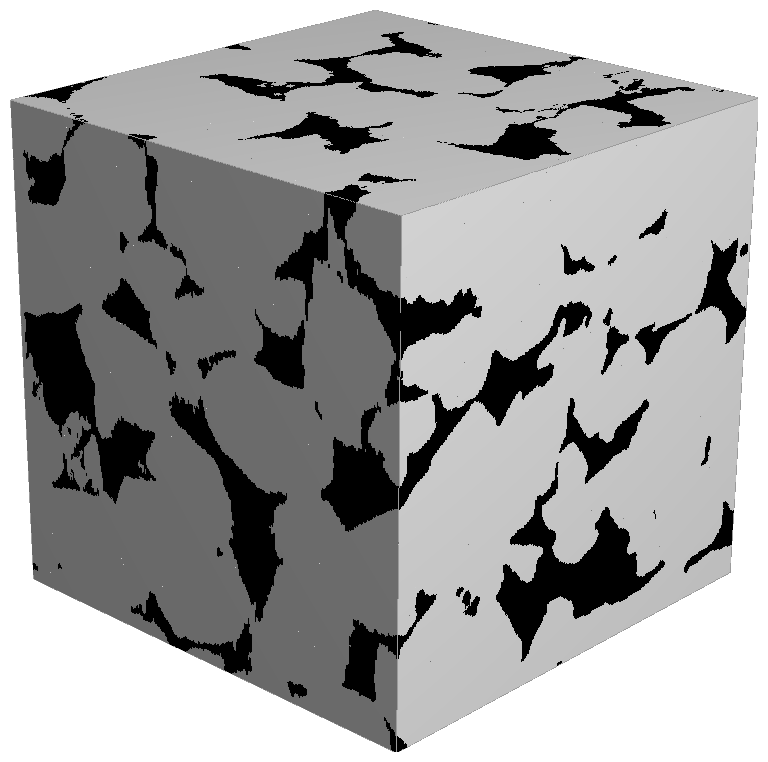}
         \caption*{ $S$}
         \label{fig:large_S}
     \end{subfigure}
     \hfill
     \begin{subfigure}[b]{0.25\textwidth}
         \centering
         \includegraphics[width=\textwidth]{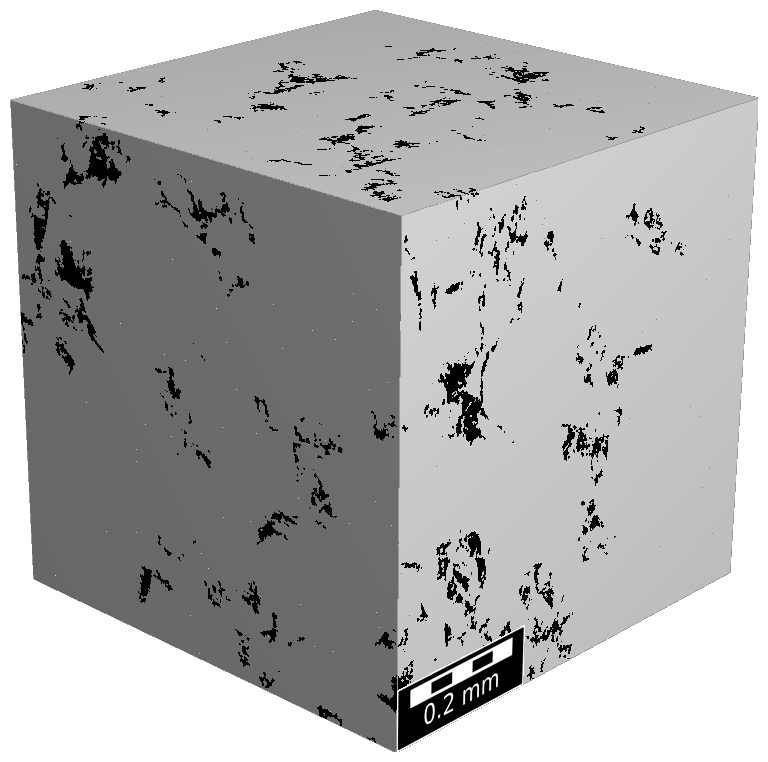}
         \caption*{$A$}
         \label{fig:large_A}
     \end{subfigure}
     \hfill
     \begin{subfigure}[b]{0.25\textwidth}
         \centering
         \includegraphics[width=\textwidth]{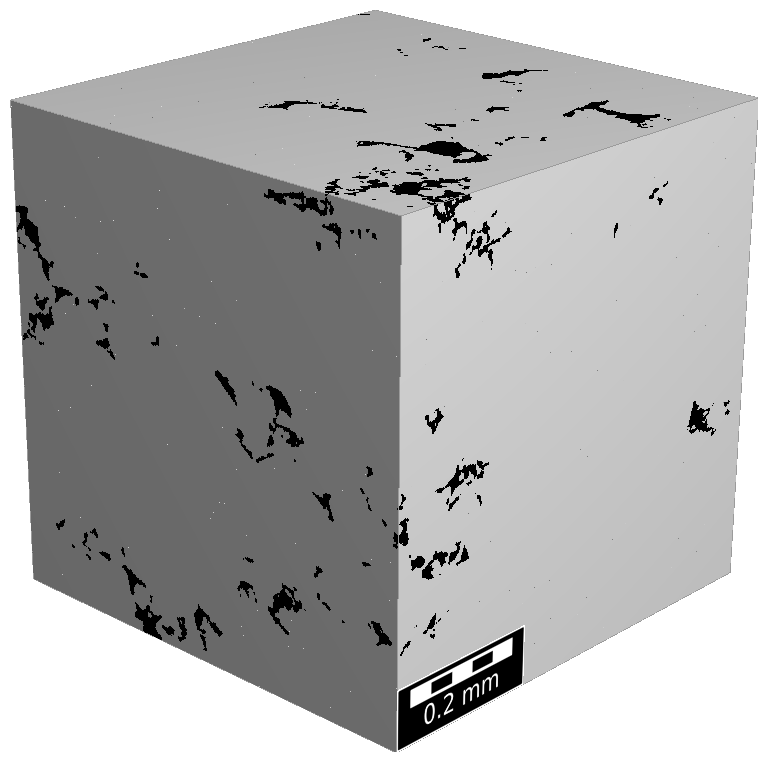}
         \caption*{ $B$}
         \label{fig:large_B}
     \end{subfigure}
     \hfill
     \begin{subfigure}[b]{0.25\textwidth}
         \centering
         \includegraphics[width=\textwidth]{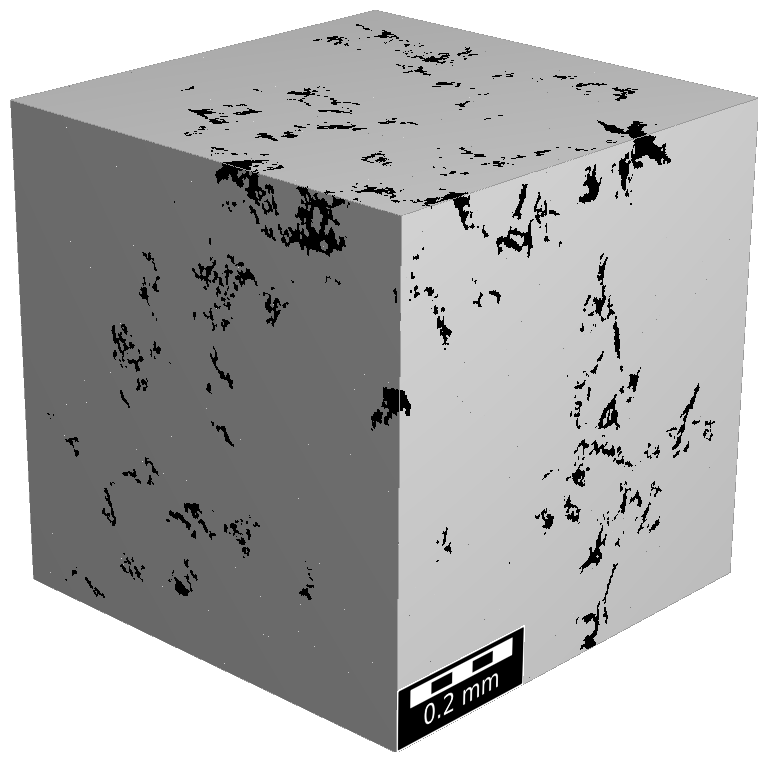}
         \caption*{ $C$}
         \label{fig:large_C}
     \end{subfigure}
     \hfill
     \begin{subfigure}[b]{0.25\textwidth}
         \centering
         \includegraphics[width=\textwidth]{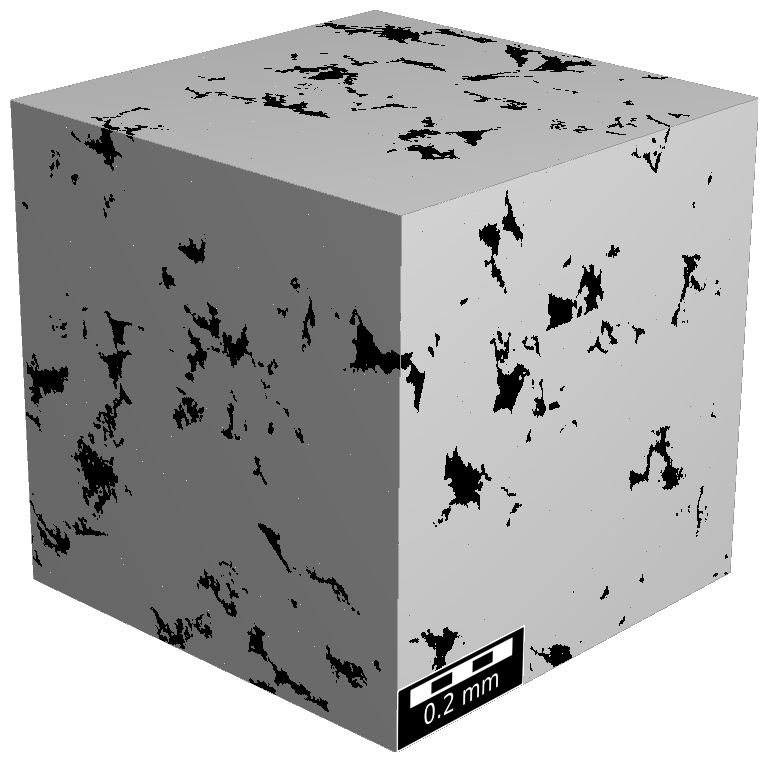}
         \caption*{ $D$}
         \label{fig:large_D}
     \end{subfigure}
     \hfill
     \begin{subfigure}[b]{0.25\textwidth}
         \centering
         \includegraphics[width=\textwidth]{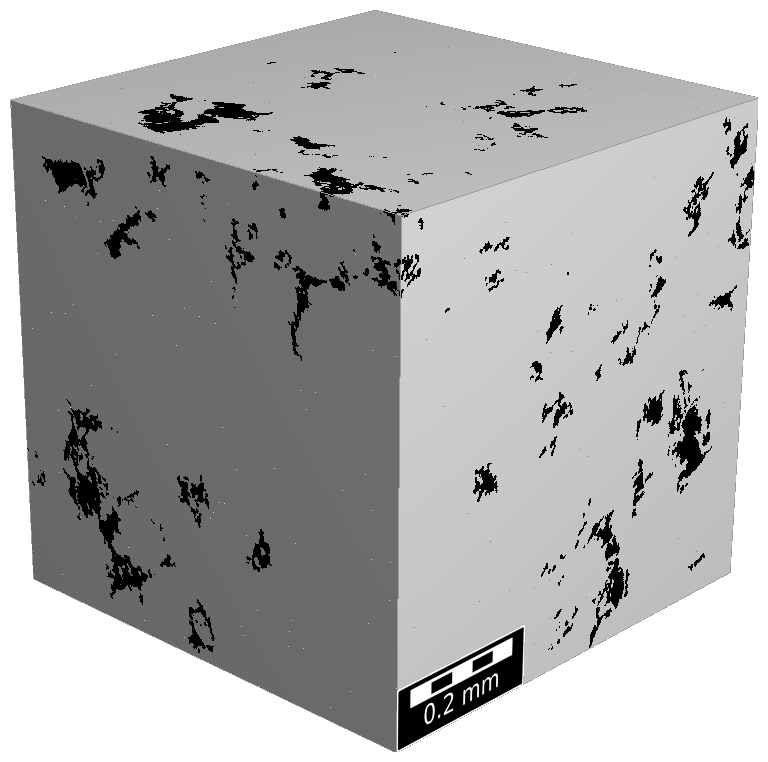}
         \caption*{ $E$}
         \label{fig:large_E}
     \end{subfigure}
        \caption{Binary images of porous media samples. The black color represents the void space $\Omega^f$, where the fluid propagates, and the grey color represents the impermeable solid. Samples A to E are ultra-tight images from \cite{orlov2021different}, and S is a moderate porosity image of Berea sandstone. For further details about the samples, see Table \ref{tab:large_samples}.}
        \label{fig:large_samples}
\end{figure}

The specific samples from our practical applications for $d=3$ are depicted in Fig. \ref{fig:large_samples}. In this context, $\Omega^f$ represents the void space of a porous medium (colored in black), and $\Gamma_{0}$ denotes the boundary between $\Omega^f$ and the solid domain, colored in grey. Along with the no-slip boundary $\Gamma_{0},$ the computational domain here formally includes exterior boundaries, i.e. the faces of the cube. However, we apply periodic boundary conditions on these faces, which do not impact the spectrum of the matrix $S$ and are therefore excluded from consideration in \eqref{stokes_bc_noslip}. The complete formulation of the flow experiment considered in our study is provided in the subsequent Section \ref{sec:problem_statement}, which includes a rigorous definition of the flow domain $\Omega^f$ and the periodic boundary conditions applied to the exterior boundaries.
It is worth mentioning , that prescribing the Dirichlet (no-slip) boundary condition on the entire boundary $\Gamma_{0}$ along with the exterior periodic conditions results in a special class of flow problems. For such flows, the pressure is determined only up to a constant nullspace, which requires special attention, as outlined later in Section \ref{sec:numerical-methods}.

If the Stokes equations \eqref{stokes} are properly discretized such that the discrete operators preserve important properties of the continuous ones, then, up to a constant in the nullspace, the Schur complement matrix $S$ is known to be spectrally equivalent to the identity operator acting on the discrete pressure space. For example, in the context of the Finite Element Method, the equivalence to the pressure mass matrix takes place with the right choice of the LBB-stable elements (see, e.g., \cite{verfurth1984combined}). 
The spectral equivalence to the identity means that the minimal nonzero eigenvalue of the Schur complement matrix $S$, as well as its effective condition number, are bounded from below and above under the mesh refinement process, i.e. when $h \to 0$. 
Moreover, in practice it is often observed that most of the eigenvalues of $S$ are equal to one (in the case of simple geometries). These facts form the basis for the famous Uzawa and Uzawa-like algorithms \cite{benzi2005numerical,elman1994inexact,bramble1997analysis,turek1999efficient,fortin2000augmented,axelsson2021krylov}, which are considered classical algorithms for solving the steady Stokes problem \eqref{stokes} and rely essentially on the preconditioning of $S$ with the identity matrix.

In fact, when solving the Stokes problem, the spectrum of $S$  depends on the boundary conditions imposed. Based on our study, we conjecture that this is due to the no-slip boundary condition \eqref{stokes_bc_noslip}, the spectrum of the Schur complement matrix contains eigenvalues which are not equal to one. Thus, when one solves problems with a small surface-to-volume ratio, which is the case most commonly considered in papers analyzing iterative solvers for Stokes problems, only a small part of eigenvalues of the Schur complement are not equal to one. This justifies using a diagonal matrix or even the identity matrix as a preconditioner.

Numerous computational studies demonstrate the efficiency of the Uzawa algorithms in the case of simple geometries. A number of reviews and theoretical studies are dedicated to this subject advancing the knowledge in the area, see, e.g. \cite{benzi2005numerical,notay2019convergence,notay2014new}.
Even a superlinear convergence of the Krylov-Uzawa algorithm can be established for general smooth geometries, see \cite{axelsson2021krylov}.
Unfortunately, for our practical application, we found that the Uzawa algorithm is not efficient.
We observe that the Schur complement matrix $S$ can become severely ill-conditioned, when computing flows in complex geometries like those representing the pore-space of tight rock sample, certain membranes, etc.. In particular, we demonstrate this issue in Section \ref{subsec:numerical_results_3d} for specific rock samples from Fig.\ref{fig:large_samples}, where the condition number of $S$ is greater than $10^5$. Therefore, development of customized methods is required for solving Stokes problems in such geometries.
The following specifics of our problem should be outlined in this context:

\begin{itemize}
    \item In our industrial application, the flow geometry $\Omega^f$ is known approximately as it is captured by Computer Tomography (CT). In this case, solving the system \eqref{coupled_matrix} with high accuracy is not necessary. Our primary objective is the computation of the absolute permeability of porous samples (see Section \ref{computing_permeability} for details), and typically, one percent accuracy in computing the objective functional is sufficient.
    \item  Methods that are efficient for solving PDEs in simple domains like squares may not perform well when applied to our complex domains.
    \item Due to limited scanning accuracy, grid convergence is usually not studied when solving PDEs in geometries derived from CT images. The goal is to efficiently solve the problem in a geometry with a fixed resolution that coincides with the voxelized CT image.
\end{itemize}

\noindent
The SIMPLE and the Uzawa preconditioners for the Schur complement matrix, denoted $\hat{S}_{\mathrm{simple}}$ and $\hat{S}_{\mathrm{uzawa}}$, are investigated in this paper. They can be written as follows:
    \begin{equation}\label{schur_uzawa_prec}
        \hat{S}_{\mathrm{simple}} = \mathbf{B} \hat{\mathbf{A}}_{\mathrm{simple}}^{-1} \mathbf{B}^T, \quad \hat{\mathbf{A}}_{\mathrm{simple}} = \mathrm{diag}(\mathbf{A}), \ \ \hat{S}_{\mathrm{uzawa}} = I.
    \end{equation}
In fact, the preconditioner $\hat{S}_{\mathrm{simple}}$ is widely-known in the CFD community since the same approximation is used in the SIMPLE  iterative method (Semi-Implicit Method for Pressure Linked Equations), which is one of the classical methods for solving the stationary Navier-Stokes equations  \cite{patankar1983calculation,patankar2018numerical,liang2016simple,elman2008taxonomy,perot1993analysis,pernice2001multigrid}, but not for the Stokes equations.

{Recently, it was demonstrated in \cite{meier2022schur} that adding discrete diffusion to the established Uzawa preconditioner significantly reduces the number of iterations in the case of channel-dominated domains. Optimal weights for the diffusion and the identity operators in the proposed preconditioner are discussed there. Broad computational experiments are performed and discussed, varying geometries, weights, finite elements spaces. It is stated that the convergence strongly depends on the geometry complexity, with the channel width and aspect ratio considered as important parameters. In our case, unlike \cite{meier2022schur}, we omit the identity. Earlier in \cite{pimanov2022workflow} we demonstrated that the diffusion-like SIMPLE preconditioner in the case of complex geometries from tight porous media performs very well. In the present article, we investigate several important aspects which were not discussed in \cite{meier2022schur} and \cite{pimanov2022workflow}. We show that the condition number of the Schur complement matrix depends on the surface-to-volume ratio, since the number of its non-unit eigenvalues is related to the surface on which no-slip boundary conditions are prescribed. Furthermore, comparing the performance of the two preconditioners, we emphasize also the fact that the SIMPLE preconditioner allows for more accurate and robust computation of the permeability. Some other differences could be mentioned. Finite element discretization is used in \cite{meier2022schur}, and iterative method (MINRES) with block diagonal preconditioner is applied for the coupled system. We use staggered finite-difference discretization of the Stokes equations, and apply the Conjugate Gradient method for the Schur complement system. As mentioned above, we identify certain geometric characteristics of the domain which are indicators for the performance of the respective Stokes solvers. \added{Only 2D synthetic problems are considered in \cite{meier2022schur}, while we perform also simulations on 3D samples from real tight reservoirs. }}

\subsection{{Paper outline and contributions.}}

The remainder of the paper is organized as follows. Section \ref{sec:problem_statement} is dedicated to the problem statement. The description of the considered iterative methods, namely the CG-Uzawa and CG-SIMPLE algorithms, is provided in Section \ref{sec:numerical-methods}. In Section \ref{sec:numerical_results}, we present and discuss the results of the computational experiments. These can be summarized as follows.

\begin{enumerate}
    \item In Section \ref{subsec:numerical_results_3d}, we investigate 3D binary samples from Fig. \ref{fig:large_samples}. Specifically, the samples with high surface-to-volume ratio coming from tight reservoirs.
    For these 3D samples, we compare the performance of the CG-SIMPLE and CG-Uzawa iterative methods and confirm by our numerical experiments that:
    \begin{itemize}
        \item  The preconditioner $\hat{S}_{\mathrm{simple}}$ provides orders of magnitude lower condition numbers than $\hat{S}_{\mathrm{uzawa}}$ hence ensuring robust and fast convergence of the CG-SIMPLE method while the CG-Uzawa method tends to stagnate;
        \item Furthermore, we demonstrate that the CG-SIMPLE provides more accurate practical computation of the absolute permeability.
    \end{itemize} 
    \item We explain this behavior in Section \ref{subsec:numerical_results_2d}, by performing a systematic study using  synthesized 2D geometries - random packings of squares. For the considered synthetic geometries, we numerically demonstrate that:
    \begin{itemize}
        \item The condition number $\mathrm{cond}(S)$ of the Schur complement matrix increases { linearly} with increasing the surface-to-volume ratio;
        \item The condition number $\mathrm{cond}(\hat{S}_{\mathrm{simple}}^{-1}S)$ of the Schur complement matrix preconditioned with the SIMPLE preconditioner decreases { super-linearly} with increasing the surface-to-volume ratio.
    \end{itemize}
    \item Additionally, for the considered synthetic geometries, we compute the full spectrum of the Schur complement matrix and observe that the number of its non-unit eigenvalues is determined by the number of boundary nodes where the Dirichlet b.c. on the tangential velocity is imposed, and by the connectivity of the flow domain. 
\end{enumerate}
Finally, conclusions are drawn.

\section{Problem statement}\label{sec:problem_statement}

In Digital Rock Physics, a typical task involves computing the absolute permeability of porous media. To compute permeability, we need to solve the steady-state Stokes problem at pore-scale resolution \cite{hornung1997homogenization,griebel2010homogenization}. The porous structure of a rock sample is typically captured using imaging techniques, such as Computer Tomography (CT).

In this section, we describe how the fluid region of a rock sample or a sample of other porous material, denoted as $\Omega^f$ in \eqref{stokes}, is represented using voxel grids. Additionally, we formulate the periodic boundary conditions typical for the DRP flow experiment and provide the formula for computing permeability.

\subsection{Representing pore-space geometries.}\label{subsec:representing_geometries}

In 3D, we represent CT images by the cubic domain $\overline{\Omega}_h = [0, L]^3$, where $L$ is the physical size of a sample. The image is decomposed into $\mathrm{n}^3$ voxels, where $\mathrm{n}$ is the number of voxels (image resolution) in each dimension:
\begin{equation}\label{cube_domain}
\overline{\Omega}_h = \bigcup_{(i,j,k) \in \mathbb{I}^{\mathrm{n}}} {\omega}_{(i,j,k)},
\end{equation}
where $\mathbb{I}^{\mathrm{n}} = \{(i,j,k): i,j,k \in {1, \ldots, \mathrm{n}}\}$ denotes a three-dimensional index set, and the voxels ${\omega}_{(i,j,k)}$ are defined as cubic regions of the length $h = L/\mathrm{n}$:
\begin{equation*}
    {\omega}_{(i,j,k)} = [(i - 1)h; ih] \times [(j - 1)h; jh] \times [(k - 1)h; kh].
\end{equation*}
The entire image $\Omega_h$ is subdivided into two parts:
\begin{equation}\label{domain_partition}
    \overline{\Omega}_h = \overline{\Omega}^f_h \cup \overline{\Omega}^s_h,
\end{equation}
that corresponds to a disjoint decomposition of the index set $\mathbb{I}^{\mathrm{n}}$:
\begin{equation}
    \mathbb{I}^{\mathrm{n}} = \mathbb{I}^{\mathrm{n}}_f \sqcup \mathbb{I}^{\mathrm{n}}_s,
\end{equation}
such that: 
\begin{equation}\label{domains_partition}
    \overline{\Omega}_h^{s} = \bigcup_{(i,j,k) \in \mathbb{I}^{\mathrm{n}}_{s}} {\omega}_{(i,j,k)}, \quad \overline{\Omega}_h^{f} = \bigcup_{(i,j,k) \in \mathbb{I}^{\mathrm{n}}_{f}} {\omega}_{(i,j,k)}.
\end{equation}
It should be noted, that arbitrary complex porous geometries can be approximated in such a way.
Example of voxel-based geometries for the case $d=3$ and $d=2$ are shown in Figs.  \ref{fig:large_samples} and \ref{fig:sol}, respectively.
In the solid region, denoted $\Omega^s_h$, the fluid does not propagate since {impermeable solid is considered here}. So, the Stokes problem \eqref{stokes} is formulated in the fluid region $\Omega^f_h$, and the no-slip boundary $\Gamma_0$ from \eqref{stokes_bc_noslip} is given as follows:
\begin{equation}\label{interior_boundary}
    \Gamma_0 = \overline{\Omega}^f_h \cap \overline{\Omega}^s_h.
\end{equation}
 It should be emphazized, that the voxel-based geometry as defined in \eqref{cube_domain} serves as the computation grid to discretize the Stokes problem in $\Omega_h^f$. The subscript $h$ here reflects the fact that the computational domain is inherently discretized since it comes from a binary CT image. In our work, we utilize the classical fully-staggered finite difference method, also known as the MAC scheme. Namely, the pressure $p^h$ is discretized at the centers of the voxels $\omega_{(i,j,k)}$, while the velocity $\mathbf{u}_h$ is discretized on the voxel faces. For a detailed description, see, for instance, \cite{harlow1965numerical,griebel1998numerical,lebedev1964difference}.

\subsection{Exterior periodic boundary conditions.}\label{subsec:bc}
In fact, apart from the interior boundary $\Gamma_0$ \eqref{interior_boundary} between solid and fluid regions, the entire boundary of $\Omega^f_h$ comprises an additional exterior part, denoted $\Gamma_{\mathrm{ext}}$.
According to the homogenization theory, periodic boundary conditions are typically imposed on the exterior faces of the domain (see, e.g., \cite{hornung1997homogenization}). As it was previously mentioned, they do not influence the spectrum of the matrix $S$.
Considering the right-handed Cartesian coordinate system, we define six boundaries, which are the outer faces of the entire cube domain $\Omega^h$ from \eqref{cube_domain}: 
  \begin{equation}
     \begin{split}
         \Gamma_{\mathrm{x}=0},  
         \quad \Gamma_{\mathrm{y}=0},
         \quad  \Gamma_{\mathrm{z}=0}, \\
        \Gamma_{\mathrm{x}=L}, 
         \quad \Gamma_{\mathrm{y}=L}, 
         \quad  \Gamma_{\mathrm{z}=L}.
     \end{split}
 \end{equation}
 Then, the exterior boundary $\Gamma_{\mathrm{ext}}$ is given as follows:
 \begin{equation}\label{boundary_ext}
     \Gamma_{\mathrm{ext}} = (\Gamma_{\mathrm{x}=0} \cup \Gamma_{\mathrm{x}=L} \cup \Gamma_{\mathrm{y}=0} \cup \Gamma_{\mathrm{y}=L} \cup \Gamma_{\mathrm{z}=0} \cup \Gamma_{\mathrm{z}=L}) \cap \overline{\Omega}_h^f.
 \end{equation}
In the DRP flow experiment, typically one direction is selected as the {\it flow direction} and two other are selected as the {\it tangential directions}.
Without loss of generality, we assume {that} $\mathrm{z}$ direction is fixed as the flow direction, and $\mathrm{x},\mathrm{y}$ directions as the tangential ones.
In the tangential directions, we impose the periodic boundary condition for both the velocity and pressure. For the $\mathrm{x}$ direction, we have:
\begin{equation*}\label{periodic}
    \mathbf{u}|_{\Gamma_{\mathrm{x}=0}} = \mathbf{u}|_{\Gamma_{\mathrm{x}=L}}, \quad 
    {p}|_{\Gamma_{\mathrm{x}=0}} = {p}|_{\Gamma_{\mathrm{x}=L}} \text{ on } \Gamma_{\mathrm{ext}},
\end{equation*}
and, similarly, the periodic boundary condition is imposed for the $\mathrm{y}$ direction.
In $\mathrm{z}$ direction, additionally a pressure difference should be accounted which drives the flow:
\begin{equation}\label{periodic_inflow}
    \mathbf{u}|_{\Gamma_{\mathrm{z}=0}} = \mathbf{u}|_{\Gamma_{\mathrm{z}=L}}, \quad 
    {p}|_{\Gamma_{\mathrm{z}=0}} = {p}|_{\Gamma_{\mathrm{z}=1}} + dp  \text{ on } \Gamma_{\mathrm{ext}},
\end{equation}
where $dp$ is a given pressure jump.
In the case of periodic boundary conditions, the geometry is assumed to be periodic too. Namely, periodicity of the flow domain $\Omega^h_f$ in all three directions can be expressed as follows:
\begin{equation}\label{periodic_geometry}
(i,j,k) \in \mathbb{I}_f^{\mathrm{n}} \implies (\mathrm{n} - i, j, k), (i, \mathrm{n} - j, k), (i, j, \mathrm{n} - k) \in \mathbb{I}_f^{\mathrm{n}}.
\end{equation}
In the case when the geometry $\mathbb{I}_f^{\mathrm{n}}$ is not periodic, it can be periodized by symmetric reflection.
It is important to emphasize, that the flow domain $\Omega_h^f$ must be connected to ensure the well-posedness of the formulated BVP \eqref{stokes}. Additionally, the constraint $\int_{\Omega_h^f} p = 0$ on the pressure is typically assumed for uniqueness of the solution.

\subsection{Computing permeability.}\label{computing_permeability}
The computation of the permeability tensor according to the homogenization theory can be found, for example in \cite{hornung1997homogenization}. {However, for brevity} we consider an approach most often used in {the} engineering literature. For the selected flow direction, $\mathrm{z}$,  the respective component of the permeability tensor, denoted $ k_{\mathrm{zz}} [m ^ 2]$,  is determined according to the Darcy's law, which is given as follows:
\begin{equation}\label{darcy_law}
    \dfrac{Q}{A} = -\dfrac{{k}_{\mathrm{zz}}}{\mu} \cdot \dfrac{dp}{L},
\end{equation}
where $\mu [Pa \cdot s]$ is the viscosity, $Q [m ^ 3 / s]$ is the volumetric flow rate, $ A [m ^ 2] $ is the cross-sectional area of the sample, $L [m]$ is the thickness of the sample in the flow direction where the pressure drop $dp[Pa]$ is applied. Note, we consider $\mu=1 [Pa \cdot s]$ since it does not influence the permeability.
The flow over the unit area is computed as follows:
\begin{equation}
    \dfrac{Q}{A} \approx \langle u_\mathrm{z} \rangle.
\end{equation}
Here, the Darcy's velocity $\langle u_\mathrm{z} \rangle$ is calculated by averaging the respective velocity component over the entire volume of the porous sample \cite{whitaker1986flow}:
\begin{equation}
    \langle u_\mathrm{z} \rangle = \frac{1}{|{\Omega}_h|}\int_{\Omega_h^f}u_\mathrm{z}.
\end{equation}
In practice, for reasons of numerical stability in finite precision arithmetic, we compute the non-dimensional permeability, denoted $\hat{k}_{\mathrm{zz}}$, by taking $L = 1 [m]$ in \eqref{cube_domain}. The physical permeability can then be computed by scaling as follows:
\begin{equation*}\label{nondim_perm}
    k_{\mathrm{zz}} = \hat{k}_{\mathrm{zz}} L^2.
\end{equation*}
When the periodic boundary conditions \eqref{periodic_inflow} are imposed in the flow direction $\mathrm{z}$, the problem is actually solved for the periodic part of the pressure, while its gradient $\nabla{p} = \frac{dp}{L}$ goes to the volumetric force $\mathbf{f}$ in \eqref{stokes}. For $\Omega_h$ having the unit length $L=1$, the unit pressure drop $dp = 1$ corresponds to the unit volume  force $\mathbf{f} = (0,0,1)^T$ applied in the flow direction \cite{wiegmann2007computation}. So, the permeability $k_{\mathrm{zz}}$ equals the Darcy (averaged) velocity $\langle u_\mathrm{z} \rangle$ in this case.
Note, that in order to compute permeabilities ${k}_{\mathrm{xx}}$ and ${k}_{\mathrm{yy}}$ for the $\mathrm{x}$ and $\mathrm{y}$ flow directions, two additional computations for $\mathbf{f} = (1,0,0)^T$ and $\mathbf{f} = (0,1,0)^T$ are required. {We omit here the discussion on computing the off diagonal elements of the permeability tensor.}

\section{Iterative methods and preconditioning}\label{sec:numerical-methods}

\subsection{Outer iterations}

Iterative methods for solving the Stokes problem \eqref{stokes} can generally be classified into two categories:
\begin{itemize}
    \item Methods in which iterations are performed for the coupled system \eqref{coupled_matrix} in its full form.
    \item Methods in which iterations are performed for the reduced Schur complement system \eqref{schur_system}.
\end{itemize}
In both categories, efficient techniques are required to solve systems involving the matrices $\mathbf{A}$ and $S$. Concerning the velocity Laplacian matrix $\mathbf{A}$, its inverse can be effectively applied through multigrid methods, such as the Algebraic Multigrid method (AMG) \cite{ruge1987algebraic}, which can handle intricate geometries. In this paper, our focus is on the Schur complement matrix $S$, so we consider the reduced formulation \eqref{schur_system}. Discussion of block preconditioners for the coupled formulation \eqref{coupled_matrix} is beyond the scope of this paper.

Since the matrix $S$ is positive semi-definite, following Axelsson (Section 3.1 in \cite{axelsson2021krylov}), we employ the Conjugate Gradient (CG) method as a Krylov subspace accelerator. However, in contrast to \cite{axelsson2021krylov}, we use the preconditioned version of CG (as described in, e.g., \cite{saad2003iterative}, Section 9.2). The specific form of the Preconditioned Conjugate Gradient (PCG) method considered in our study is formulated in {Algorithm} \ref{alg_pcg}.
\begin{algorithm}
\caption{Preconditioned CG, adapted from Axelsson \cite{axelsson2021krylov}}
\label{alg_pcg}
\begin{algorithmic}[1]
\Require tolerance $\varepsilon_S$, initial guess $p_h^0$
\Ensure Approximate solution $p_h$  for the system \eqref{schur_system}
\State Compute the initial residual $r_h^0 = Sp_h^0 - g_h$
\State Solve $\hat{S}z_h^0 = r_h^0$ for $z_h^0$
\State Set $d_h^0 = z_h^0$, $k = 0$
\While{not converged}
    \State Apply $Sd_h^k = q_h^k$ for $q_h^k$
    \State $\alpha_k = \dfrac{(r_h^k)^T z_h^k}{(d_h^k)^T q_h^k}$
    \State $p_h^{k+1} = p_h^k + \alpha_k d_h^k$
    \State $r_h^{k+1} = r_h^k + \alpha_k q_h^k$
    \State Solve $\hat{S}z_h^{k+1} = r_h^{k+1}$ for $z_h^{k+1}$
    \State  If $\|z_h^{k+1}\| / \|z_h^{0}\| < \varepsilon_S$, exit loop
    \State $\beta_k = \dfrac{(r_h^{k+1})^T z_h^{k+1}}{(r_h^k)^T z_h^k}$
    \State $d_h^{k+1} = z_h^{k+1} + \beta_k d_h^k$, $k = k + 1$
\EndWhile
\State Return $p_h^{k+1}$ as the approximate solution
\end{algorithmic}
\end{algorithm}
By \textbf{{CG-Uzawa}} and \textbf{{CG-SIMPLE}}, we denote the {Algorithm} \ref{alg_pcg} for the preconditioners $\hat{S} = \hat{S}_{\mathrm{uzawa}}$ and $\hat{S} = \hat{S}_{\mathrm{simple}}$ defined in \eqref{schur_uzawa_prec}, respectively. Obviously, when the identity is used as a preconditioner, the preconditioned and unpreconditioned CG coincide. It should be noted, that in the Stokes case, the spectrum of $\hat{S}_{\mathrm{simple}}$ is qualitatively different from the spectrum of $S$. Namely, the matrix $\hat{S}_{\mathrm{simple}}$ behaves essentially as the pressure Laplacian matrix $\mathbf{B}\mathbf{B}^T$, so its condition number increases quadratically as the grid resolution decreases.
However, such spectral behavior turns out to be justified in the case of tight geometries in the presence of narrow channels where the nature of flow is predominantly diffusive.

\begin{remark}
   In the original work by Axelsson, the CG-Uzawa algorithm is formulated for the regularized version of the system \eqref{coupled_matrix}. We do not employ any regularization as there is no necessity for it.
\end{remark}

\begin{remark}
    As mentioned earlier, in the case of periodic boundary conditions \eqref{periodic_inflow}, the matrix $S$ has a constant nullspace. However, the CG method converge to the normal solution as soon as the corresponding system is consistent \cite{kaasschieter1988preconditioned}.
    Because $S$ is singular, instead of considering its inverse, we have to actually consider its pseudo-inverse. For clarity of notation, we use $\hat{S}^{-1}$ instead of $\hat{S}^{\dagger}$ throughout the text. The term condition number here is used in the sense of effective condition number.
\end{remark}

\noindent
According to {Algorirhm} \ref{alg_pcg}, the preconditioned residual norm is used as the stopping criteria for the outer iterative process, see line 10 there. Specifically, given an input tolerance $\varepsilon_{S}$, the outer PCG iterations stop as soon as:
\begin{equation}\label{unpreconditioned_residual}
    \|\hat{S}^{-1}r_h^{\#}\| / \|\hat{S}^{-1}r_h^0\| < \varepsilon_S,
\end{equation}
where the superscript $\#$ denotes the final iteration number, and $r_h^k$ denotes the residual on the $k^{th}$ outer iteration, given as follows:
\begin{equation}\label{residual_outer}
    r_h^k = Sp_h^k - g_h.
\end{equation}
Using the preconditioned residual norm is considered more natural because the preconditioned CG method minimizes the norm of the preconditioned residuals at each step \cite{saad2003iterative,van2003iterative}. However, in our numerical experiments, for the CG-SIMPLE method, we additionally monitor the unpreconditioned residual norm.
This helps us compare with the CG-Uzawa, where the unpreconditioned norm is monitored. 
Moreover, when studying synthetic 2D geometries in Section \ref{subsec:numerical_results_2d}, we use the unpreconditioned residual norm directly as the stopping criteria  for consistency of the comparison; Note, {similar approach was used in \cite{meier2022schur}.}

\subsection{Inner iterations and stopping criteria}

On each step of the outer iterations, applying the matrix $S$ requires solution of an auxiliary problem to recover intermediate velocity from the intermediate pressure. We use the PCG method for solving with the velocity Laplacian matrix $\mathbf{A}$, so formally we deal with inexact version of the outer PCG method \cite{bramble1997analysis,elman1994inexact}. Thus, we have a two-level inner-outer iterative process: at each step of the outer PCG iteration for $S$, inner PCG iterations for $\mathbf{A}$ are performed.
In our numerical experiments, we use preconditioned relative residual norm as the stopping criteria for the inner iterations with matrix $\mathbf{A}$. Namely, given an input tolerance $\varepsilon_{\mathbf{A}}$, the inner PCG iteration for computing $\mathbf{u}_h = \mathbf{A}^{-1}\mathbf{f}^h$ stops as soon as:
\begin{equation}
        \|(\hat{\mathbf{A}})^{-1}r_\mathbf{A}^{\#}\| / \|(\hat{\mathbf{A}})^{-1} r_\mathbf{A}^0\| < \varepsilon_{\mathbf{A}},
\end{equation}
where $r_\mathbf{A}^{k}$ denotes the residual on the $k^{th}$ inner iteration, given as follows:
\begin{equation}
    r_{\mathbf{A}}^k = \mathbf{A} \mathbf{u}_h^k - \mathbf{f}_h.
\end{equation}
Additionally, at each step of the  outer CG-SIMPLE iteration, we have to solve the system with the preconditioner matrix $\hat{S}_{\mathrm{simple}}$. Again, we use the Preconditioned Conjugate Gradient for solving with the preconditioner $\hat{S}$. It should be noted that, as well as the pressure Schur complement matrix $S$, the pressure Neumann Laplacian  $\hat{S}_{\mathrm{simple}}$ has constant in the nullspace, so a special care is required. As for the stopping criteria, we use preconditioned relative residual norm  for inner iterations with the matrix $\hat{S}_{\mathrm{simple}}$. For example, given an input tolerance $\varepsilon_{\hat{S}}$, the inner CG iteration for computing $p_h = (\hat{S}_{\mathrm{simple}})^{-1}g^h$ stops as soon as:
\begin{equation}
        \|(\hat{\hat{S}}_{\mathrm{simple}})^{-1} r_{\hat{S}}^{\#}\| / \|(\hat{\hat{S}}_{\mathrm{simple}})^{-1} r_{\hat{S}}^0\| < \varepsilon_{\hat{S}},
\end{equation}
where $r_{\hat{S}}^{k}$ denotes the residual on the $k^{th}$ inner iteration, given as follows:
\begin{equation}
    r_{\hat{S}}^{k} = (\hat{S}_{\mathrm{simple}}) p_h^k - g_h.
\end{equation}
For building preconditioners $\hat{\mathbf{A}}$ and $\hat{\hat{S}}_{\mathrm{simple}}$ for the Laplacian matrices $\mathbf{A}$ and $\hat{S}_{\mathrm{simple}}$, we use the implementation BoomerAMG \cite{yang2002boomeramg} from HYPRE library. 


According to Algorithm \ref{alg_pcg}, we do not employ a flexible version of the Krylov subspace method for the outer iterative process, although we do not compute exact solutions during the inner iterations in our practical applications. This approach aligns with the consensus among numerous researchers, who assert that if the inner iterations are solved with sufficient accuracy, there is no necessity for a flexible version of the outer iterative method.
Moreover, we consider Axelsson's assertion \cite{axelsson2021krylov}, Section 3.1, p.615, which claims that the CG-accelerated inner-outer Uzawa algorithm, used for solving the Schur complement formulation of the Stokes problem, is not sensitive to the tolerance used for inner iterations. At least it is less sensitive then the classical, not accelerated Uzawa algorithm.
Our observations also confirm that the CG-SIMPLE algorithm converges well when a reasonable tolerance for the inner iterations is maintained. 

\section{Computational experiments}\label{sec:numerical_results}

Computational experiments are conducted to numerically investigate:
\begin{itemize}
    \item Possible correlation between the surface-to-volume ratio and the condition number of the {preconditioned\textbackslash unpreconditioned} Schur complement matrix;
    \item Possible correlation between the number of boundary nodes where no-slip boundary conditions are imposed and the number of non-unit eigenvalues of the Schur complement matrix.
    \item The performance of the Uzawa and SIMPLE preconditioners in solving the Schur complement problem \eqref{schur_system}, especially in complex geometries with high surface-to-volume ratio;
    \item The performance of two preconditioners in computing the permeability of representative 3D and 2D samples according to \eqref{darcy_law};
\end{itemize}




\subsection{3D rock samples from tight reservoirs}\label{subsec:numerical_results_3d}

\subsubsection{\it Preliminaries: samples and notations} 

In this section, we study the performance of the CG-Uzawa and CG-SIMPLE algorithms for the pore space images of real rock samples from tight reservoirs. We consider six 3D images: five ultra-tight samples $A-E$ earlier considered in \cite{orlov2021different}, and one image of the classical Berea's sandstone with medium porosity \cite{pimanov2022workflow}. The corresponding pore space images are depicted on Fig. \ref{fig:large_samples} using the Geodict visualization tool. The samples $A-E$ are scanned with resolution $1.2$ mkm and have the size $\mathrm{n} = 600$, and the sample S is scanned with the resolution $4$ mkm and has the size $\mathrm{n} = 300$. Detailed information about the samples can be found in Table \ref{tab:large_samples}, which includes the reference permeability ${k}_{\mathrm{zz}}^{\mathrm{ref}}$ computed according to \eqref{darcy_law} solving \eqref{schur_system} for the respective samples with the CG-SIMPLE with very high accuracy, the number of fluid voxels $\mathbb{V}^f = |\mathbb{I}^{\mathrm{n}}_f|$, and the porosity $\phi$, defined as follow:

\begin{equation}\label{porosity}
    \phi = (\mathbb{V}^f / \mathbb{V}) \cdot 100 \%,
\end{equation}
where $\mathbb{V} = |\mathbb{I}^{\mathrm{n}}| = \mathrm{n}^d$ is the total number of voxels. 
Additionally, for each sample we compute the surface-to-volume ratio which is defined as follows:
\begin{equation}\label{s-t-v-ration}
    \sigma_s = (\mathbb{V}_{surf}^s / \mathbb{V}^f) \cdot 100\%,
\end{equation}
where the surface area $\mathbb{V}_{surf}^s$ of the no-slip boundary is determined as the number of near-boundary solid voxels, i.e., solid voxels face-adjacent with the fluid domain:
\begin{equation}\label{surface_area}
    \mathbb{V}_{surf}^s = |\mathbb{I}^{\mathrm{n}}_{surf}|, \quad \mathbb{I}^{\mathrm{n}}_{surf} = \{(i,j,k) \in \mathbb{I}^{\mathrm{n}}_s: \omega_{(i,j,k)} \cap \overline{\Omega}_f \neq \emptyset\}.
\end{equation}
Note, the surface area $\mathbb{V}_{surf}^s$ coincides with the number of boundary nodes (i.e.,  the nodes lying on $\Gamma_0$), where the Dirichlet b.c. on the tangential velocity component is imposed.

Recall that in \cite{pimanov2022workflow}, the computation of the permeability for the samples $A-E$ with the current CG-SIMPLE algorithm was compared to computations with other commercial and academic codes, and it was shown that all the values are in the same range.

\begin{table}[h!]
    \centering
    \caption{Description of 3D samples $A-S$ including the problem size $\mathrm{n}$, reference permeabilities ${k}_{\mathrm{zz}}^{\mathrm{ref}}$, the number of fluid voxels $\mathbb{V}^f$, the porosity $\phi$, and the surface-to-volume ratio $\sigma_s$.}
    \begin{tabular}{||c | c c c c c : c||} 
     \hline
      & \textbf{A} & \textbf{B} & \textbf{C} & \textbf{D} & \textbf{E} & \textbf{S} \\
     \hline
     size $\mathrm{n}$ & 600 & 600 & 600 & 600 & 600 & 300 \\
     perm. ${k}_{\mathrm{zz}}^{\mathrm{ref}}$, mD & $0.65$ & $0.78$ & $0.31$ & $9.6$ & $0.75$ & $6.61 \cdot 10^{3}$ \\
    \# pores $\mathbb{V}^f$, mln. & $12.5$ & $9.4$ & $11.4$ & $20.9$ & $13.9$ & $5.7$ \\
     porosity $\phi$, \% & $5.8$ & $4.4$ & $5.3$ & $9.7$ & $6.4$ &  $21.1$ \\
     s-t-v ratio $\sigma_s$, \% & 74 & 56  & 70 & 55 & 61 & 28   \\
     \hline
    \end{tabular}
    \label{tab:large_samples}
\end{table}

\subsubsection{\it Performance of the preconditioners in solving the Schur complenment problem and in computing the permeability}\label{sec:3d_performance}

As the stopping criteria for the outer CG iterations, we use $\varepsilon_S = 10^{-3}$ for both the CG-Uzawa and CG-SIMPLE algorithms. The relative {preconditioned} residual norm, as defined in \eqref{unpreconditioned_residual}, is used here.   For the inner iterations, we use a higher precision $\varepsilon_{\mathbf{A}} = 10^{-6}$ in both cases. Also, for the CG-SIMPLE algorithm we use $\varepsilon_{\hat{S}} = 10^{-6}$ for solving with the SIMPLE preconditioner $\hat{S}_{\mathrm{simple}}$. The reference permeabilities $k_{\mathrm{zz}}^{\mathrm{ref}}$ were computed using the CG-SIMPLE algorithm with higher precision $\varepsilon_S = 10^{-5}$, $\varepsilon_{\mathbf{A}} = 10^{-8}$, $\varepsilon_{\hat{S}} = 10^{-8}$. Such inexact solves for $\mathbf{A}$ and $\hat{S}_{\mathrm{simple}}$ make it difficult to rigorously analyze the underlying numerical methods. However, the purpose of this section is to demonstrate convergence problems with the established Uzawa preconditioner that occur in practical permeability calculations when computing flows in tight porous media. So, we compare methods for the settings that we usually use in our practical calculations. 
\begin{figure}
     \centering
     \caption{Convergence history of the CG-SIMPLE and CG-Uzawa algorithms for 3D samples $A-S$ described in the Table \ref{tab:large_samples}. Preconditioned relative residual norm $\|\hat{S}^{-1}r_S^{k}\| / \|\hat{S}^{-1}r_S^0\|$ (top) and the permeability error $e^k$ defined in \eqref{scaled_permeability_error} (bottom) are shown.}
      \begin{subfigure}[b]{\textwidth}
         \centering
         \includegraphics[width=\textwidth]{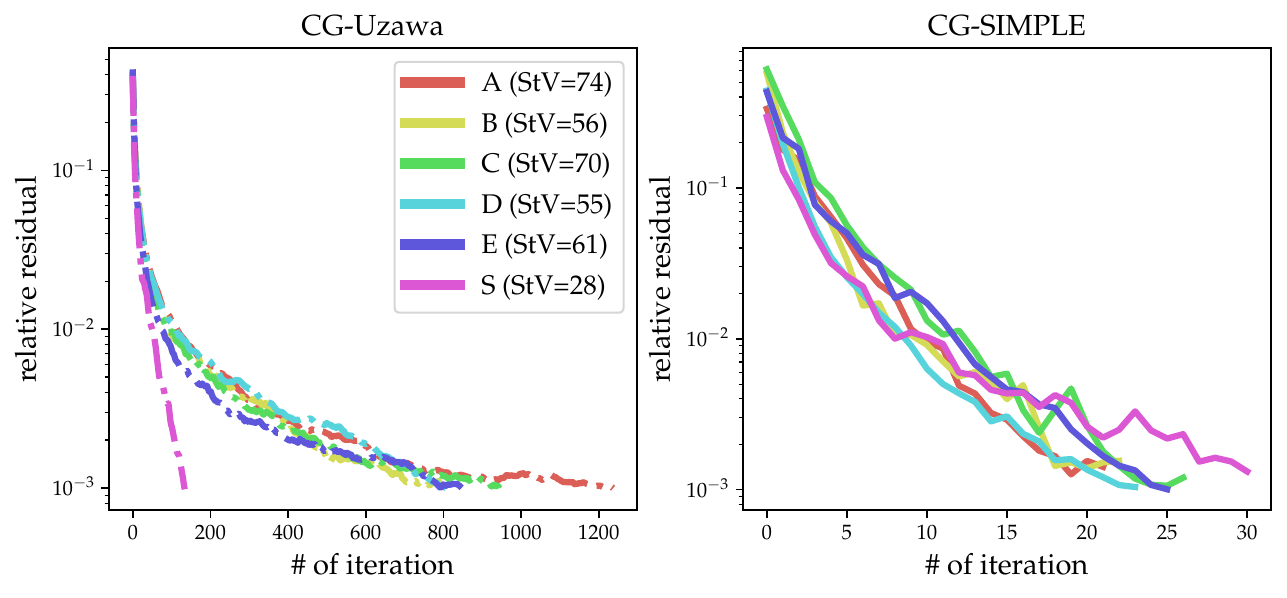}
     \end{subfigure}
     
     \begin{subfigure}[b]{\textwidth}
         \centering
            \includegraphics[width=\textwidth]{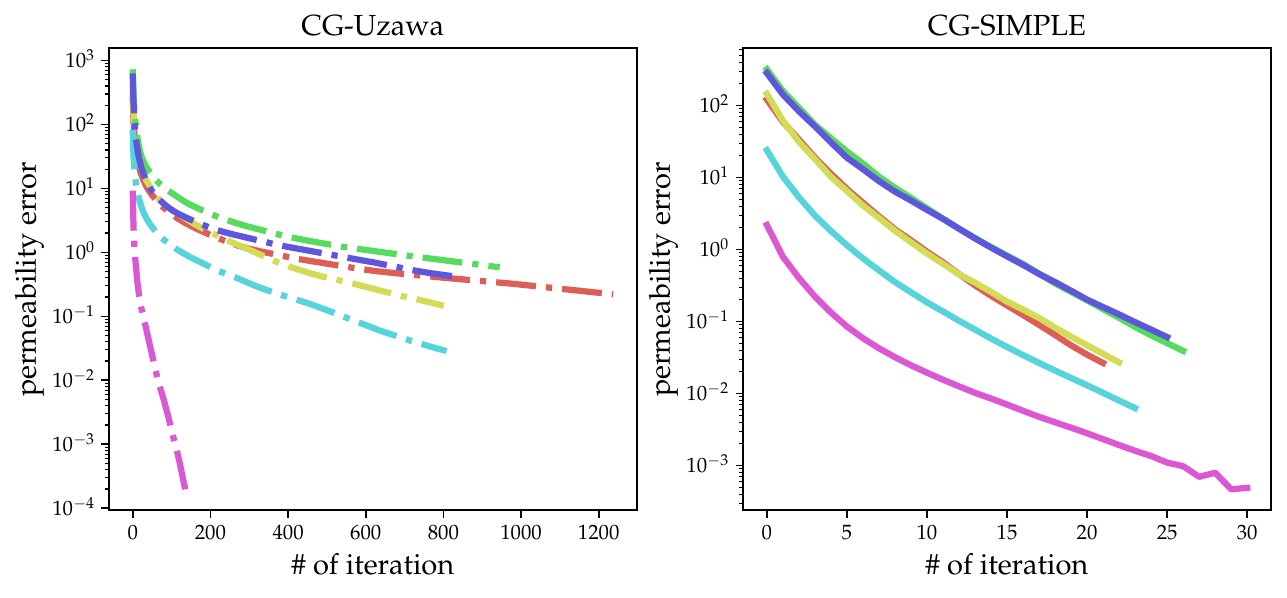}
     \end{subfigure}
    \label{fig:large_res_perm}
\end{figure}
The convergence history for the selected tolerances is presented in Fig. \ref{fig:large_res_perm}, where the relative unpreconditioned residual norm $\|{\hat{S}^{-1}}r_S^{k}\| / \|{\hat{S}^{-1}}r_S^0\|$ is shown on the top graph.
Furthermore, the relative permeability error, defined as:
\begin{equation}\label{scaled_permeability_error}
    e^k = |{k_{\mathrm{zz}}^k} - k_{\mathrm{zz}}^{\mathrm{ref}}|  / {k_{\mathrm{zz}}^{\mathrm{ref}}},
\end{equation}
is shown in Fig. \ref{fig:large_res_perm} on the bottom graph.
The convergence summary can be found in Table \ref{tab:large_sample_convergence}, which includes the number of iterations, the required computational time, and the permeability error $e^{\#}$ computed on the final iteration. Additionally, we provide estimations for the condition numbers of the preconditioned and unpreconditoned Schur complement matrices, computed for free during the outer CG iterations using Lanczos algorithm. 
\begin{table}[h!]
    \centering
    \caption{Summary of the results for CG-Uzawa and CG-SIMPLE algorithms for 3D samples $A-S$ described in the Table \ref{tab:large_samples} including relative permeability error on the final iteration $e^{\#}$, the total computational time, the number of iterations, and the estimated condition numbers.}
\begin{tabular}{|| c | c c c c c : c||} 
 \hline
  & \textbf{A} & \textbf{B} & \textbf{C} & \textbf{D} & \textbf{E} & \textbf{S} \\
 \hline
  \textbf{CG-Uzawa}: &  &  &  &  &  &  \\
    perm. error $e^{\#}$, \%  & 0.221 & 0.146 & 0.585 & 0.029 & 0.403 & 0.0002 \\
    total comp. time, hrs & $8.5$ & $5.0$ & $6.2$ & $10.9$ & $7.3$ &  $0.8$ \\
    $\#$ iters & $1238$ & $802$ & $946$ & $809$ & $849$ & $136$ \\
     $\approx\mathrm{cond}(S)$ & $7.1 \cdot 10^{5}$  & $2.1 \cdot 10^{5}$& $3.7 \cdot 10^{5}$ & $2.2 \cdot 10^{5}$ & $3.4 \cdot 10^{5}$ & $3.4 \cdot 10^{3}$\\
 \hline\hline
 \textbf{CG-SIMPLE}: &  &  &  &  &  &  \\
    perm. error $e^{\#}$, \% & 0.019 & 0.020 & 0.031 & 0.005 & 0.047 & 0.0005 \\
    total comp. time, hrs & $0.3$ & $0.3$ & $0.3$ & $0.7$ & $0.5$ & $0.3$ \\
    $\#$ iters & $22$ & $23$ & $27$ & $24$ & $26$ & $31$ \\
    $\approx\mathrm{cond}(\hat{S}_{\mathrm{simple}}^{-1}S)$ & $0.9 \cdot 10^{2}$  & $1.2 \cdot 10^{2}$ & $1.4 \cdot 10^{2}$ & $1.1 \cdot 10^{2}$ &$1.8 \cdot 10^{2}$ & $2.4  \cdot 10^{2}$ \\
 \hline
\end{tabular}
    \label{tab:large_sample_convergence}
\end{table}
The following hardware was used in our numerical experiments: 48x MPI compute node (Dell PowerEdge M640), dual Intel Xeon Gold 6132 ("Skylake") @ 2.6 GHz, i.e. 28 CPU cores per node. The computational times shown in Table \ref{tab:large_sample_convergence} were obtained using 8 CPU nodes.

Several observations can be drawn from the results presented in Table \ref{tab:large_sample_convergence}. For the considered low porosity, high surface-to-volume ratio images, the SIMPLE preconditioner appears to perform better than the Uzawa preconditioner, as it converges robustly while the CG-Uzawa {tends to} stagnate. Firstly, despite the CG-SIMPLE algorithm being more expensive (approximately x1.5) per iteration, its total computational time is smaller compared to the CG-Uzawa because significantly fewer number of iterations is required. Secondly, the estimated condition number for the preconditioned Schur complement matrix is about three orders of magnitude smaller for the SIMPLE preconditioner than for the Uzawa preconditioner. For moderate porosity (sample $S$), the condition number for both preconditioners is comparable. Thirdly, the CG-SIMPLE computes the permeability much more accurately for the considered samples. Actually, achieving even 10\% accuracy in computing permeability is not always possible with the established CG-Uzawa method.

It is worth mentioning, that for the CG-SIMPLE algorithm oscillations may appear in the permeability error (see Fig. \ref{fig:large_res_perm}, sample $S$). This is the case when inner tolerance $\varepsilon_{\hat{S}}$ for inverting preconditioner $\hat{S}=\hat{S}_{\mathrm{simple}}$ is not small enough. However, despite these oscillations, the permeability error continues to decrease, albeit not monotonically. This confirms our conjecture that using flexible Krylov method is not necsesary in our case. Note, similar oscillations may appear if inner tolerance $\varepsilon_{\mathbf{A}}$ for inverting $\mathbf{A}$ when applying $S$ is not small enough.
It is also worth mentioning that the unpreconditioned residual for the CG-SIMPLE algorithm \textit{always} decreases monotonically in our numerical experiments (i.e. even when oscillations appear in the permeability error). In Fig. \ref{fig:large_res_perm_unprec_simple} in the Appendix, we show the unpreconditioned residual for the CG-SIMPLE algorithm corresponding to the stopping criteria $\varepsilon_{S} = 10^{-3}$ in the preconditioned residual norm, as it is shown in Fig. \ref{fig:large_res_perm}. Note, the unpreconditioned residual norm is suppressed only up to $10^{-2}$ in this case. We investigate this issue in the subsequent section for synthetic 2d geometries, where we {use} the unpreconditioned residual norm directly as the stopping criteria {for the preconditioned CG} for both methods. Detailed comparison of preconditioned/unpreconditioned residuals for CG-SIMPLE and CG-Uzawa for different thresholds $\varepsilon_{S} = 10^{-1}$, $\varepsilon_{S} = 10^{-2}$, $\varepsilon_{S} = 10^{-3}$ is shown in Tables \ref{table:01}, \ref{table:001}, \ref{table:0001} in the Appendix, respectively.

\subsubsection{\it Correlation between the surface-to-volume ratio and the condition number of the Schur complement matrix} 

\begin{figure}
     \centering
         \caption{Correlation between surface-to-volume ratio and estimated condition number for the samples $A-S$ described in the Table \ref{tab:large_samples}.}
         \includegraphics[width=0.5\textwidth]{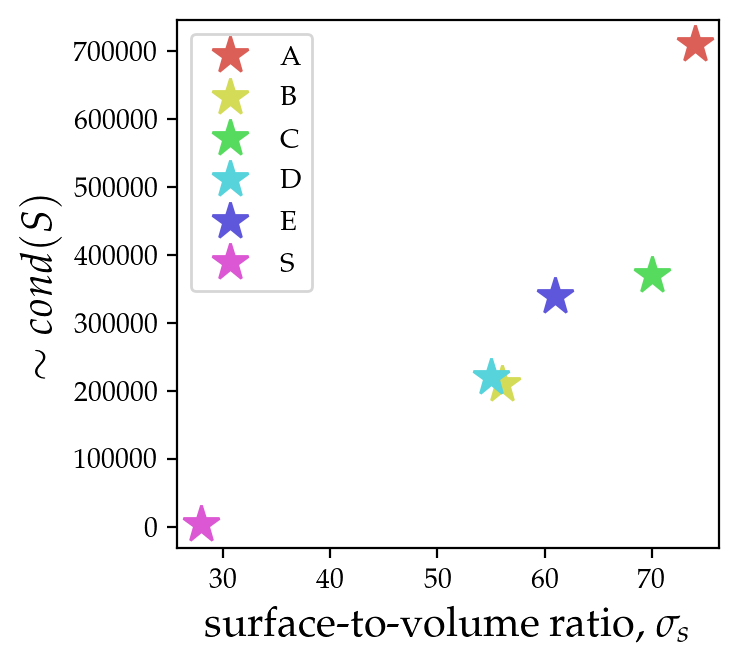}
         \label{fig:stv-cond}
\end{figure}

In Fig. \ref{fig:stv-cond}, we show that for the considered samples $A-S$ there is a strong correlation between the surface-to-volume ratio $\sigma_s$ \eqref{s-t-v-ration} and the estimated condition number of the unpreconditioned Schur complement matrix. This is also confirmed by the number of iterations for the CG-Uzawa algorithm.
However, the dependence in the case of the Schur complement matrix preconditioned with the SIMPLE is not so distinctive here, so we perform a more rigorous study in the subsequent Section \ref{sec:2d_correlation}.

\subsection{
Two-dimensional synthetic geometries}\label{subsec:numerical_results_2d}
 To identify consistent patterns that affect the performance of the methods under consideration in complex pore space domains, we consider synthetic 2D geometries with a transparent generation process and directly available geometric information (such as porosity, no-slip surface area, etc). As in the 3D case, we present and discuss the convergence of the two algorithms, as well as the accuracy with which the permeability is computed. Additionally, for the considered synthetic samples we compute and analyse the full spectra of the preconditioned and unpreconditioned Schur complement matrices. \\

\noindent
\subsubsection{Generation of synthetic 2D geometries.} 
We study flows passing around arrays of solid square obstacles randomly placed in a fluid bed. First of all, a uniform voxel (pixel in 2D) grid is generated in $\Omega$ as described in Section \ref{sec:problem_statement}. For ease of generation, we consider square obstacles of the same size, and each obstacle is located in the center of the square cell, or is slightly shifted, so that a cell contains the obstacle and a part of the flow domain around it. The obstacles do not touch the boundary of the cell.  Each generated geometry is defined by four integer parameters $(N,n_\mathrm{c},n_{\mathrm{avg}},n_{\min})$, where $N$ and $n_\mathrm{c}$ determine the number of cells in one direction and their size measured in voxels, while the parameters $n_{\mathrm{avg}}$ and $n_{\min}$ control the average and minimal thicknesses (in voxels) of the fluid channels between two adjacent solid squares. Note, that the cells and the obstacles in all cases are adjusted to the introduced computational grid, so that each voxel is fully occupied either by fluid or by solid.
\begin{figure}
     \centering
    \caption{Example of synthetic 2D geometry: array of randomly shifted square obstacles for $N = 7$, $n_\mathrm{s} = 40$, $n_{\mathrm{avg}} = 10$, $n_{\min} = 2$.}
     \begin{subfigure}[b]{0.45\textwidth}
         \centering
         \includegraphics[width=\textwidth]{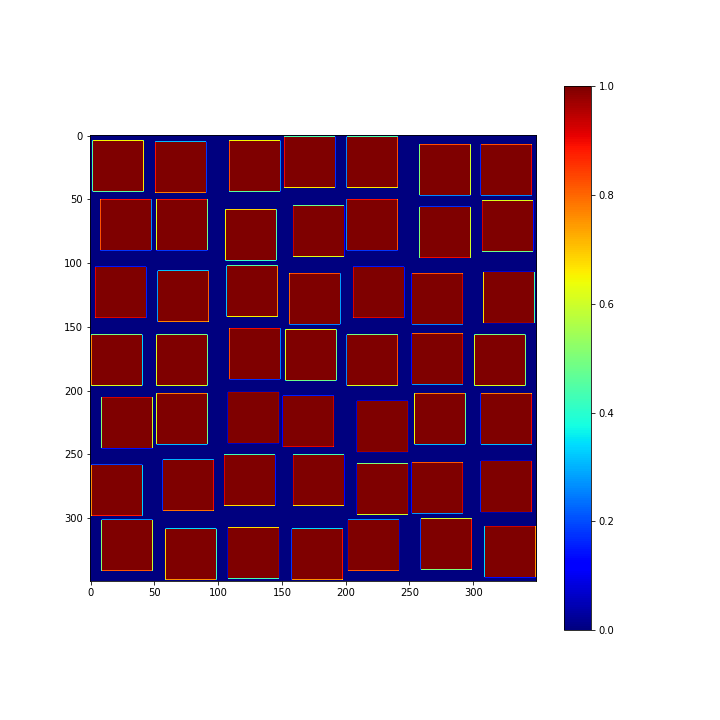}
         \caption{geometry}
         \label{fig:sol_geom}
     \end{subfigure}
     \begin{subfigure}[b]{0.44\textwidth}
         \centering
         \includegraphics[width=\textwidth]{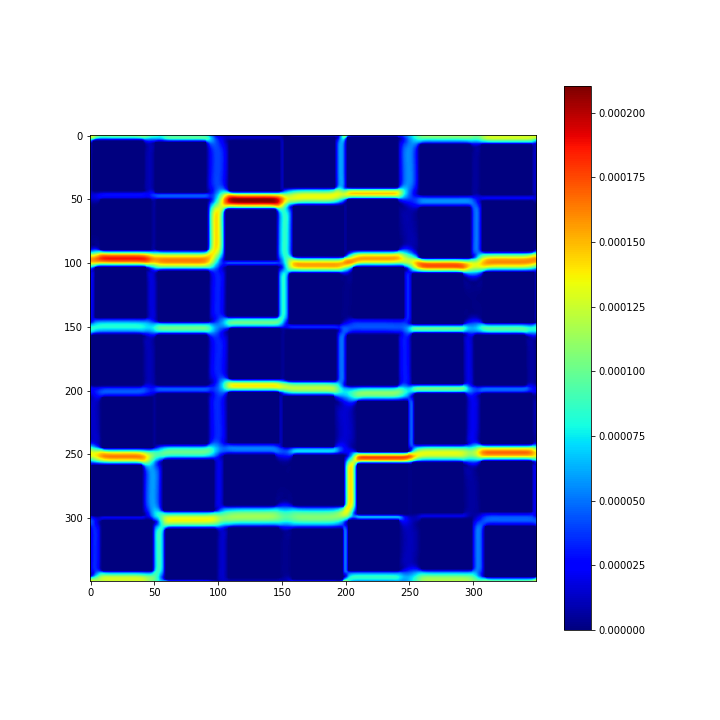}
         \caption{velocity magnitude}
         \label{fig:sol_velocity}
     \end{subfigure}
        \label{fig:sol}
\end{figure}
An example for $N = 7$, $n_{\mathrm{c}} = 50$, $n_{\mathrm{avg}} = 10$, $n_{\min} = 2$ is presented in Fig. \ref{fig:sol_geom}. Formally, we have a domain of the total size $\mathrm{n}\times \mathrm{n}$, $\mathrm{n} = Nn_\mathrm{c}$ which represents an $N \times N$ array of cells of the size $n_\mathrm{c} \times n_\mathrm{c}$; each cell contains a solid square of the size $(n_\mathrm{c} - n_\mathrm{avg}) \times (n_\mathrm{c} - n_\mathrm{avg})$ with the origin $(n_\mathrm{c}/2 + r_1, n_\mathrm{c}/2 + r_2)$, where $r_1, r_2 \in [-(n_{\mathrm{avg}}-n_{\min})/2, (n_{\mathrm{avg}}-n_{\min})/2]$ are (integer) random shifts. It should be noted, that randomness is necessary to observe non-trivial solutions which take place in the case of fully periodic arrays. 


\subsubsection{\it Performance of the preconditioners in solving the Schur complement problem and in computing  permeability} 

In the present section, we investigate the effect of surface-to-volume ratio on the convergence of the algorithms by varying {the} average thickness {of the channels} of the synthetic 2D geometries. Namely, we randomly generated five geometries according to the procedure described in the previous subsection for $N = 7$, $n_{\mathrm{c}} = 50$, $n_{\min} = 2$, and $n_{\mathrm{avg}} = \{4,6,8,10,12\}$.
Detailed information about the samples can be found in Table \ref{tab:2d_avg_description}, which includes the reference permeability ${k}_{\mathrm{xx}}^{\mathrm{ref}}$, the number of fluid voxels $\mathbb{V}^f$, the porosity $\phi$ \eqref{porosity}, and the surface-to-volume ratio $\sigma_s$ \eqref{s-t-v-ration}. It should be noted, that the average channel thickness for synthetic 2D geometries is directly related to the surface-to-volume ratio for general porous media. 

\begin{table}[h!]
    \centering
    \caption{Description of the synthetic 2D geometries with variable channel thicknesses $n_{\mathrm{avg}}$ including the problem size $\mathrm{n}$, reference permeabilities ${k}_{\mathrm{xx}}^{\mathrm{ref}}$, the number of fluid voxels $\mathbb{V}^f$, the porosity $\phi$, and the surface-to-volume ratio $\sigma_s$.}
    \begin{tabular}{||c | c c c c c||} 
\hline
      & $n_{\mathrm{avg}}=4$ & $n_{\mathrm{avg}}=6$ & $n_{\mathrm{avg}}=8$ & $n_{\mathrm{avg}}=10$ & $n_{\mathrm{avg}}=12$ \\
\hline
     size $\mathrm{n}$ & 350 & 350 & 350 & 350 & 350 \\
     perm. $k_{\mathrm{xx}}^{\mathrm{ref}}$ & $1.01 \cdot 10^{-6}$ & $3.17 \cdot 10^{-6}$ & $7.31 \cdot 10^{-6}$  & $1.50 \cdot 10^{-5}$ & $2.53 \cdot 10^{-5}$ \\
     \# pores $\mathbb{V}^f$, thsnd. & $18.8$ & $27.6$ & $36.0$ & $44.1$ & $51.8$  \\
     porosity $\phi$, \% & 15.4 & 22.6 & 29.4 & 36.0 & 42.3  \\
     s-t-v ratio $\sigma_s$, \% &  46.9 & 30.5 & 22.3 & 17.3 & 14.0   \\
     \hline
    \end{tabular}
    \label{tab:2d_avg_description}
\end{table}

\begin{figure}
     \centering
      \caption{Convergence history of the CG-SIMPLE and CG-Uzawa algorithms for the synthetic 2D geometries with variable channel thicknesses $n_{\mathrm{avg}}$ described in the Table \ref{tab:2d_avg_description}. Unpreconditioned relative residual norm $\|r_h^{k}\| / \|h_S^0\|$ (right) and relative permeability error $e^k$ defined in \eqref{scaled_permeability_error} (left) are shown. }
     \includegraphics[width=\textwidth]{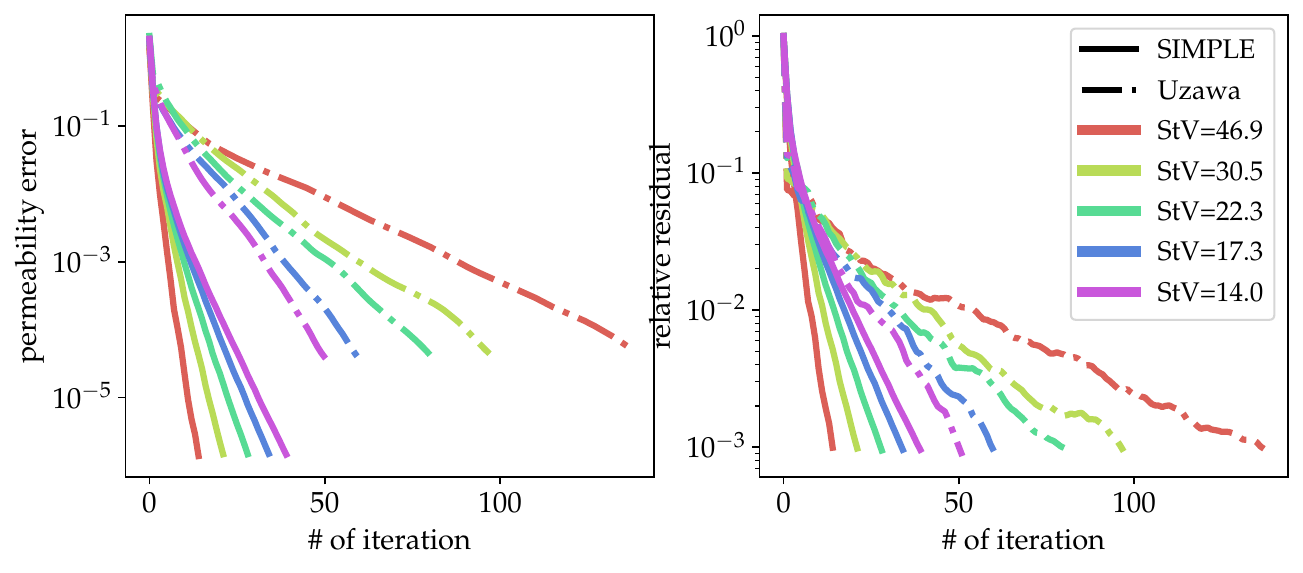}
\label{fig:2D_convergence}
\end{figure}

  As the stopping criteria for the outer CG iterations, we use  $\varepsilon_S = 10^{-3}$ for both the CG-SIMPLE and CG-Uzawa algorithms. In contrast to the previous subsection, we consider the unpreconditioned relative residual here.
  Considering residuals in the same norm allows for straightforward comparison between the two algorithms. A detailed comparison of preconditioned/unpreconditioned residuals for the CG-SIMPLE and CG-Uzawa for different thresholds $\varepsilon_{S} = 10^{-2}$, $\varepsilon_{S} = 5\cdot10^{-3}$, $\varepsilon_{S} = 10^{-3}$ is shown in Tables \ref{table:2d_01}, \ref{table:2d_005}, \ref{table:2d_001} in the Appendix, respectively. As for the inner iterations, in this experiment we use the machine epsilon $\varepsilon_{\mathbf{A}} = \varepsilon_{\hat{S}} = 10^{-13}$. In particular, this means that the exact Uzawa is used here. 
 The convergence history for selected tolerances is presented in Fig. \ref{fig:2D_convergence}, where the relative unpreconditioned residual norm $\|r_h^{k}\| / \|r_h^0\|$ is shown on the left and the permeability error $e^k$ defined in \eqref{scaled_permeability_error} is shown on the right.
The convergence summary can be found in Table \ref{tab:2d_avg_summary}, which includes the number of iterations until convergence, the required computational time, and the permeability error $e^{\#}$ computed on the final iteration. In addition, we show the effective condition numbers of the preconditioned and unpreconditoned Schur complement matrices, which are computed using a direct method. 
Similar observations to those made for the 3D simulations can be made from the results presented in Table \ref{tab:2d_avg_summary}. Again, the CG-SIMPLE algorithm requires less iterations and less total computational time compared to the CG-Uzawa for samples with larger surface-to-volume ratio. When decreasing surface-to-volume ratio, the CG-Uzawa tends to show better performance. As in the 3D case, the CG-SIMPLE more accurately computes the permeability, which is clearly seen from the left graph in Fig. \ref{fig:2D_convergence}.
This may seem counterintuitive since the residuals have the same norm for both algorithms, as illustrated in the right graph of Fig. \ref{fig:2D_convergence}. This phenomenon requires a deeper theoretical investigation, which is beyond the scope of this paper. We currently hypothesize that the primary reason is that the iterative solutions computed by the two algorithms belong to different Krylov subspaces. Furthermore, it is evident from Tables \ref{table:2d_01}, \ref{table:2d_005}, \ref{table:2d_001} that as the tolerance $\varepsilon_S$ for the outer iterations decreases, the difference in permeability calculations decreases as well.
\begin{table}[h!]
    \centering
    \caption{Convergence summary of the CG-Uzawa and CG-SIMPLE algorithms for the synthetic 2D geometries with variable channel thicknesses $n_{\mathrm{avg}}$ described in the Table \ref{tab:2d_avg_description} including the permeability error on the final iteration $e^{\#}$, the total computational time, the number of iterations, and the condition numbers.}
\begin{tabular}{|| c | c c c c c ||} 
\hline
      & $n_{\mathrm{avg}}=4$ & $n_{\mathrm{avg}}=6$ & $n_{\mathrm{avg}}=8$ & $n_{\mathrm{avg}}=10$ & $n_{\mathrm{avg}}=12$ \\
  & (StV=$46.9$) & (StV=$30.5$) & (StV=$22.3$) & (StV=$17.3$) & (StV=$14.0$)  \\
 \hline
  \textbf{CG-Uzawa}: &  &  &  &  &   \\
    perm. error $e^{\#}$, \%  & $5.4 \cdot 10^{-5}$ & $4.4 \cdot 10^{-5}$ & $4.3 \cdot 10^{-5}$ & $3.6 \cdot 10^{-5}$ & $3.4 \cdot 10^{-5}$  \\
    total comp. time, s & $6.5$ & $6.0$ & $5.8$ & $5.0$ & $4.7$  \\
    $\#$ iters & $138$ & $98$ & $81$ & $61$ & $52$  \\
     $\mathrm{cond}(S)$ & $4.7 \cdot 10^{3}$  & $2.3 \cdot 10^{3}$& $1.3 \cdot 10^{3}$ & $7.7 \cdot 10^{2}$ & $5.3 \cdot 10^{2}$ \\
 \hline\hline
 \textbf{CG-SIMPLE}: &  &  &  &  &   \\
    perm. error $e^{\#}$, \% & $1.4 \cdot 10^{-6}$ & $1.5 \cdot 10^{-6}$ & $1.5 \cdot 10^{-6}$ & $1.5 \cdot 10^{-6}$ & $1.5 \cdot 10^{-6}$  \\
    total comp. time, s & $1.4$ & $2.4$ & $3.5$ & $4.6$ & $5.7$  \\
    $\#$ iters & $15$ & $22$ & $29$ & $35$ & $40$  \\
    $\mathrm{cond}(\hat{S}_{\mathrm{simple}}^{-1}S)$ & $3.4 \cdot 10^{1}$  & $8.6 \cdot 10^{1}$ & $1.6 \cdot 10^{2}$ & $2.6 \cdot 10^{2}$ &$3.4 \cdot 10^{2}$  \\
 \hline
\end{tabular}
    \label{tab:2d_avg_summary}
\end{table}
It is important to note that for the 2D geometries under consideration, the convergence of the CG-Uzawa algorithm closely follows a linear asymptote. Specifically, the convergence factor is determined by the effective condition number:
\begin{equation}\label{cond_cg_rate}
    \dfrac{\sqrt{\mathrm{cond}(S)} - 1}{\sqrt{\mathrm{cond}(S)} + 1}.
\end{equation}
Taking into account $\lambda_{\max}(S) = 1$, \eqref{cond_cg_rate} becomes:
\begin{equation}\label{lambda_cg_rate}
    \dfrac{1-\sqrt{\lambda_{\min}(S)}}{1+\sqrt{\lambda_{\min}(S)}},
\end{equation}
where $\lambda_{\min}(S) > 0$ is the smallest non-zero eigenvalue.
Notably, for the CG-SIMPLE, the convergence is super-linear, which means some clustering of eigenvalues. {Similar spectral clustering is reported in \cite{meier2022schur}.}

\subsubsection{\it Correlation between the surface-to-volume ratio and the condition number of the Schur complement matrix}\label{sec:2d_correlation}
\begin{figure}
     \centering
      \caption{Correlation between surface-to-volume ratio and condition number for the synthetic 2D samples with various channel thicknesses $n_{\mathrm{avg}}$ described in the Table \ref{tab:2d_avg_description}.}
     \begin{subfigure}[b]{0.45\textwidth}
         \centering
         \includegraphics[width=\textwidth]{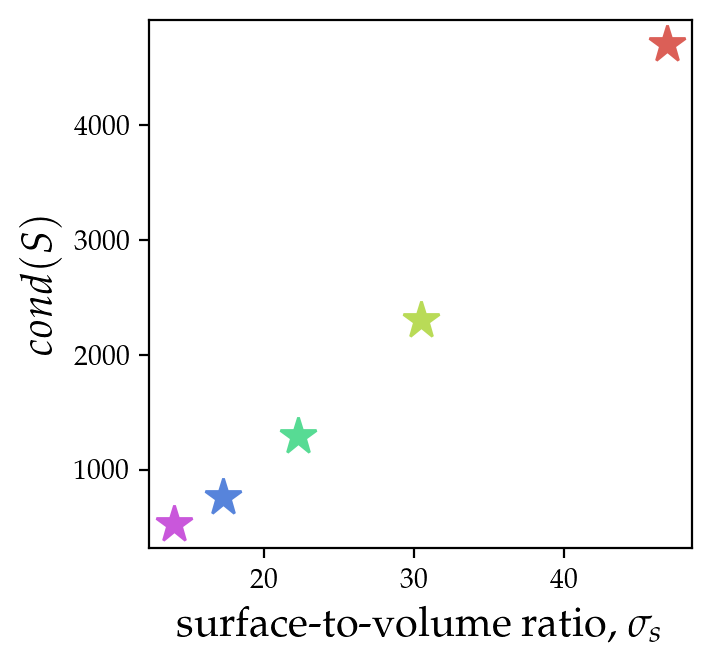}
         \caption{Uzawa}
     \end{subfigure}
          \begin{subfigure}[b]{0.45\textwidth}
         \centering
         \includegraphics[width=\textwidth]{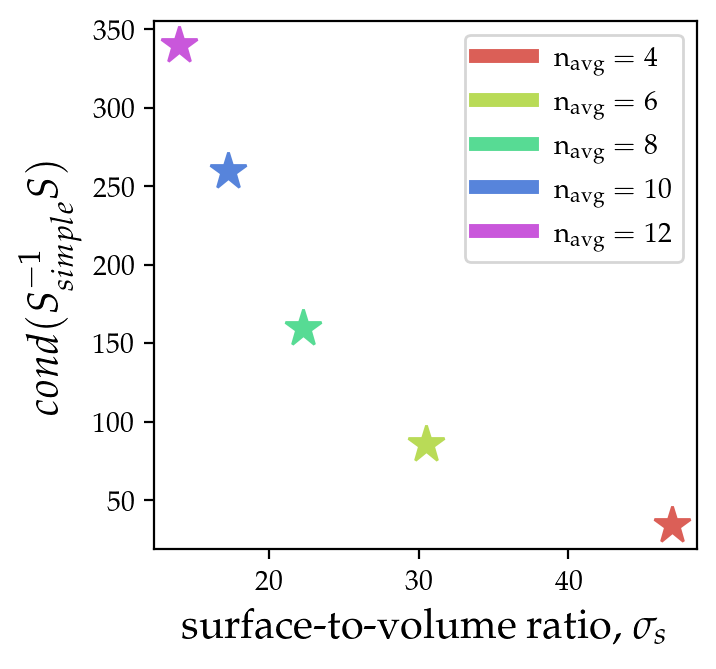}
         \caption{SIMPLE}
     \end{subfigure}
\label{fig:2d_stv_cond}
\end{figure}

In Fig. \ref{fig:2d_stv_cond}, for the considered 2D geometries with variable channel thicknesses, we show the dependence between the surface-to-volume ratio $\sigma_s$ \eqref{s-t-v-ration} and the condition number of the unpreconditioned/preconditioned Schur complement matrix.
The main conclusion that we draw is that the condition number of the Schur complement matrix increases linearly when increasing the surface-to-volume ratio.
 However, unlike the 3D case, for the considered here synthetic 2D geometries, the inverse dependence is also clearly observed for the Schur complement matrix preconditioned with the SIMPLE. Specifically, the condition number of the preconditioned matrix decreases super-linearly
 with increasing surface-to-volume ratio.
Moreover, we conjecture that the condition number of the Schur complement matrix is influenced by the number of non-unit eigenvalues in its spectrum. In the next subsection, we show that the large ratio of non-unit eigenvalues in the spectrum of the Schur complement matrix is directly related to the large surface-to-volume ratio. 

\subsubsection{\it Number of non-unit eigenvalues and no-slip surface area.}

\begin{figure}
     \centering
         \caption{Dependence of the number of non-unit eigenvalues  of the Schur complement matrix $N_{ev}$ on the connectivity of the flow domain and the surface area of the no-slip boundary $\mathbb{V}_{surf}^s$ defined in \eqref{surface_area}. The number of square obstacles $N$ in one dimension is shown by color, total $N^2$ obstacles.}
         \includegraphics[width=0.5\textwidth]{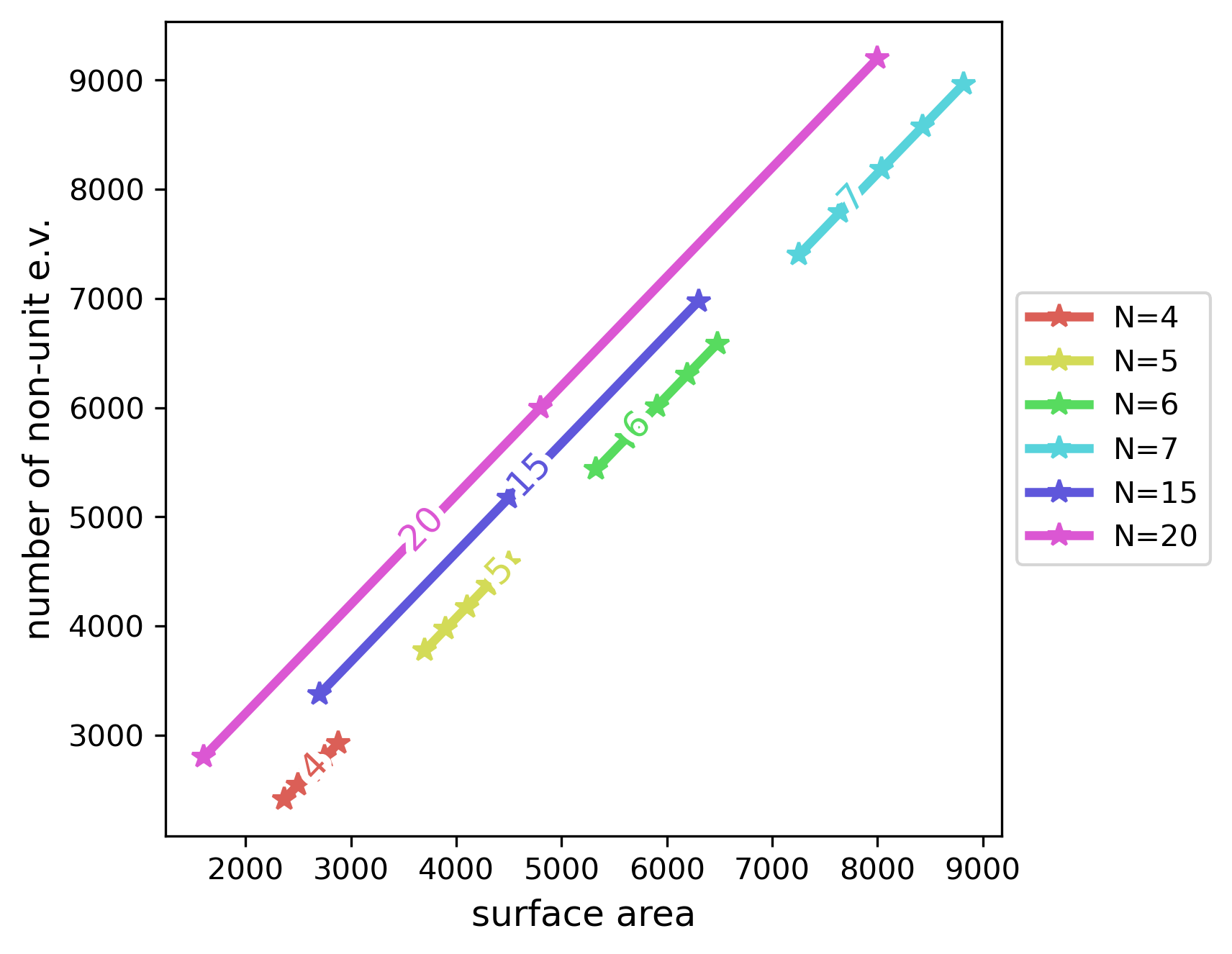}
         \label{fig:spectrum}
\end{figure}

In the present section, for the synthetic 2D geometries, we show the relation between surface-to-volume ratio and the number of non-unit eigenvalues of the Schur complement matrix.
We consider various configurations of synthetic 2D geometries. First, we vary the number of squares $N$, which determines the degree of connectivity of the flow domain $\Omega_h^f$. Second, for each $N$ we vary the surface area $\mathbb{V}_{surf}^s$ defined in \eqref{surface_area}.
For the selected configurations, we compute the full spectrum of the Schur complement matrix $S$ and calculate the number of its eigenvalues, denoted $N_{ev}$, which are not equal to one (including the zero eigenvalue). The results are presented in Figure \ref{fig:spectrum}, where we show the dependence between the surface area $\mathbb{V}_{surf}^s$ and the number of non-unit eigenvalues $N_{ev}$.
For the considered here class of geometries, we observe that the following empirical formula for the number of non-unit eigenvalues holds:
\begin{equation}\label{spectrum_relation}
    N_{ev} = \mathbb{V}_{surf}^s + 3N^2 - 1.
\end{equation}
In particular, the formula \eqref{spectrum_relation} reveals that fot a simply-connected domain, the surface-to-volume ratio equals the ratio of non-unit to unit eigenvalues of the Schur complement matrix. Thus, the greater is the ratio, the further the Schur complement is from the identity, and the worse is the performance of the Uzawa preconditioner.

\section{Conclusions}


In conclusion, the article presents a comparative study of the CG-SIMPLE and CG-Uzawa algorithms for solving the Stokes equations in tight geometries{, together with explanation of the reasons for good or bed performance of the later}. The results show that the CG-SIMPLE algorithm applied for 3D samples of rocks from tight reservoirs  ensures robust and fast convergence to high accuracy both in solving the Schur complement problem and in computing the permeability, while the established CG-Uzawa algorithm tends to stagnate. This behavior is further explained through a systematic study of  synthetic 2D geometries. The main conclusion that we draw is that for the considered here synthetic geometries, the condition number of the Schur complement matrix increases linearly with increasing the surface-to-volume ratio.
 Meanwhile, the condition number of the Schur complement matrix preconditioned with the SIMPLE preconditioner decreases super-linearly with increasing surface-to-volume ratio. We also demonstrate that the number of non-unit eigenvalues of the Schur complement matrix is determined by the number of boundary nodes where the Dirichlet b.c. on the tangential velocity component is imposed, and by the connectivity of the flow domain. These findings provide important insights into the behavior of the solvers for the Schur complement matrix and demonstrate effectiveness of the SIMPLE preconditioner for solving the Stokes problem in tight geometries.

\section*{Acknowledgements}
Ivan Oseledets was supported by Alexander von Humboldt Research Award. Oleg Iliev was supported by BMBF under grant 05M20AMD ML-MORE.

\bibliographystyle{plain}      
\bibliography{main.bib}   

\newpage
\appendix

\subsection{Convergence in preconditioned and unprecinditioned norm.}

The Appendix provides additional information for numerical experiments from Sections \ref{sec:3d_performance} and \ref{sec:2d_correlation}. Tables \ref{table:01}, \ref{table:001}, \ref{table:0001} correspond to the convergence plot presented in Fig. \ref{fig:large_res_perm}. The tables include relative permeability error $e^k$ and the corresponding number of iteration $k$ for different residual thresholds $\varepsilon_{S}=10^{-1}$, $\varepsilon_{S}=10^{-2}$, $\varepsilon_{S}=10^{-3}$, respectively. Note, for the CG-SIMPLE algorithm, the tables include thresholds for both preconditioned and unpreconditioned residual norms for the comparison purposes. Similar results corresponding to the convergence plot presented in Fig. \ref{fig:2D_convergence} for thresholds $\varepsilon_{S}=10^{-2}$, $\varepsilon_{S}=5\cdot 10^{-3}$, $\varepsilon_{S}=10^{-3}$ are provided in Tables \ref{table:2d_01}, \ref{table:2d_005}, and \ref{table:2d_001}, respectively.

\begin{figure}
     \centering
     \caption{Unpreconditioned residual norm for the CG-SIMPLE algorithm corresponding to the convergence history from Fig. \ref{fig:large_res_perm}.}
    \includegraphics[width=0.5\textwidth]{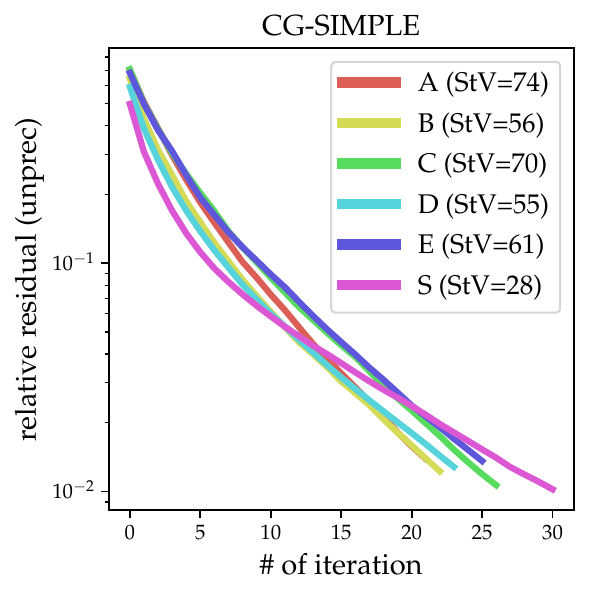}
    \label{fig:large_res_perm_unprec_simple}
\end{figure}

Additionally, in Fig. \ref{fig:large_res_perm_unprec_simple}, we plot the unpreconditioned residual norm for the CG-SIMPLE algorithm corresponding to the convergence history from Fig. \ref{fig:large_res_perm}, where the preconditioned residual is plotted. Note, the unpreconditioned residual always decreases monotonically.

\begin{table}[h!]
\centering
 \begin{tabular}{|c | c | c | c|} 
 \hline
   sample & Uzawa & SIMPLE (prec) & SIMPLE (unprec) \\
 \hline
 A & 33.24 (9) & 16.86 (3) & 1.829 (8) \\
 B & 51.87 (8) & 19.28 (3) & 1.955 (8)\\
 C & 76.97 (9) & 39.10 (4) & 4.161(10)\\
 D & 11.20 (8) & 3.320 (3) & 0.571 (7)\\
 E & 53.79 (8) & 51.42 (3) & 3.565(10)\\
 S & 0.7281(6) & 0.4033(2) & 0.059 (6)\\
 \hline
\end{tabular}
\caption{Permeability error $e^k$ (iteration $k$) for residual threshold $\varepsilon_S = 0.1$ corresponding to the convergence history from Fig. \ref{fig:large_res_perm}.}
\label{table:01}
\end{table}

\begin{table}[h!]
\centering
 \begin{tabular}{|c | c | c | c|} 
 \hline
   sample & Uzawa & SIMPLE (prec) & SIMPLE (unprec) \\
 \hline
 A & 3.412(114) & 0.840 (10) & 0.0128 (23) \\
 B & 5.013 (92) & 0.678 (11) & 0.0162 (24) \\
 C & 9.430 (88) & 1.613 (13) & 0.0357 (27) \\
 D & 1.119(112) & 0.278 (9) & 0.0044 (25) \\
 E & 6.097 (74) & 1.427 (13) & 0.0960(29) \\
 S & 0.022 (51) & 0.016 (11) & 0.0005 (31) \\
 \hline
\end{tabular}
\caption{Permeability error $e^k$ (iteration $k$) for residual threshold $\varepsilon_S = 0.01$ corresponding to the convergence history from Fig. \ref{fig:large_res_perm}.}
\label{table:001}
\end{table}

\begin{table}[h!]
\centering
 \begin{tabular}{|c | c | c | c|} 
 \hline
   sample & Uzawa & SIMPLE (prec) & SIMPLE (unprec) \\
 \hline
 A & 0.2209(1237) & 0.0169 (22) & 0.00120 (39) \\
 B & 0.0146 (800) & 0.0221 (23) & 0.00260 (44) \\
 C & 0.5852 (944) & 0.0557 (25) & 0.00260 (48) \\
 D & 0.0286 (807) & 0.0084 (22) & 0.00006 (47) \\
 E & 0.4030 (837) & 0.0936 (26) & 0.01230 (55) \\
 S & 0.0002 (134) & 0.0005 (31) & 0.00000 (58) \\
 \hline
\end{tabular}
\caption{Permeability error $e^k$ (iteration $k$) for residual threshold $\varepsilon_S = 0.001$ corresponding to the convergence history from Fig. \ref{fig:large_res_perm}.}
\label{table:0001}
\end{table}

\begin{table}[h!]
\centering
 \begin{tabular}{|c | c | c | c |} 
 \hline
    & Uzawa & SIMPLE (prec) & SIMPLE (unprec) \\
 \hline
    s-t-v=46.9 & $6.6e-3$ (55) & $1.1e-4$ (8) & $1.1e-4$ (8)  \\
    s-t-v=30.5 & $4.7e-3$ (42) & $4.8e-4$ (9) & $1.1e-4$ (12) \\
    s-t-v=22.3 & $5.6e-3$ (33) & $7.7e-4$ (11) & $1.4e-4$ (15) \\
    s-t-v=17.3 & $3.8e-3$ (30) & $11.7e-4$ (11) & $1.5e-4$ (18) \\
    s-t-v=14.0 & $3.8e-3$ (25) & $9.5e-4$ (13) & $1.6e-4$ (20) \\
 \hline
\end{tabular}
\caption{Permeability error $e^k$ (iteration $k$) for residual threshold $\varepsilon_S = 0.01$ corresponding to the convergence history from Fig. \ref{fig:2D_convergence}.}
\label{table:2d_01}
\end{table}

\begin{table}[h!]
\centering
 \begin{tabular}{|c | c | c | c |} 
 \hline
    & Uzawa & SIMPLE (prec) & SIMPLE (unprec)  \\
 \hline
    s-t-v=46.9 & $2.1e-3$ (76)& $1.1e-4$ (8) & $2.24e-5$ (10)  \\
    s-t-v=30.5 & $1.7e-3$ (53)& $1.6e-4$ (11) & $2.67e-5$ (15)  \\
    s-t-v=22.3 & $1.4e-3$ (47)& $2.5e-4$ (14) & $3.17e-5$ (19)  \\
    s-t-v=17.3 & $1.1e-3$ (38)& $1.9e-4$ (17) & $3.25e-5$ (23)  \\
    s-t-v=14.0 & $0.7e-3$ (35)& $2.5e-4$ (18) & $3.64e-5$ (26)  \\
 \hline
\end{tabular}
\caption{Permeability error $e^k$ (iteration $k$) for residual threshold $\varepsilon_S = 0.005$ corresponding to the convergence history from Fig. \ref{fig:2D_convergence}.}
\label{table:2d_005}
\end{table}

\begin{table}[h!]
\centering
 \begin{tabular}{|c | c | c | c |} 
 \hline
    & Uzawa & SIMPLE (prec) & SIMPLE (unprec)  \\
 \hline
    s-t-v=46.9 & $5.39e-5$ (137) & $9.6e-6$ (12) & $2.24e-6$ (14)  \\
    s-t-v=30.5 & $4.43e-5$ (97) & $5.9e-6$ (17) & $1.47e-6$ (21)  \\
    s-t-v=22.3 & $4.26e-5$ (80) & $1.4e-5$ (22) & $1.47e-6$ (28)  \\
    s-t-v=17.3 & $3.62e-5$ (60) & $1.24e-5$ (27) & $1.49e-6$ (34) \\
    s-t-v=14.0 & $3.40e-5$ (51) & $1.43e-5$ (29) & $1.46e-6$ (39) \\
 \hline
\end{tabular}
\caption{Permeability error $e^k$ (iteration $k$) for residual threshold $\varepsilon_S = 0.001$ corresponding to the convergence history from Fig. \ref{fig:2D_convergence}.}
\label{table:2d_001}
\end{table}

\end{document}